\documentclass{opt2025} 
\usepackage{float}
\usepackage{amsmath, amssymb}
\usepackage{booktabs}
\usepackage{ragged2e}
\usepackage{array}
\newcolumntype{Y}{>{\raggedright\arraybackslash}X}

\newcommand{\Real}{\mathbb{R}}

\newcommand{\E}{\mathbb{E}}

\newcommand{\prox}{\operatorname{prox}}
\newcommand{\hjprox}{\operatorname{prox}^{\delta}}

\newcommand{\amp}{\mathop{\:\:\,}\nolimits}

\title[HJ-Prox-based Operator Splitting]{A Monte Carlo Approach for Nonsmooth Convex Optimization via Proximal Splitting Algorithms}

\optauthor{%
\Name{Nicholas Di} \Email{nd56@rice.edu}\\
\addr Rice University
\AND
\Name{Eric C. Chi} \Email{echi@umn.edu}\\
\addr University of Minnesota
\AND
\Name{Samy Wu$\text{ }$Fung} \Email{swufung@mines.edu}\\
\addr Colorado School of Mines
}

\begin{document}

\maketitle

\begin{abstract}%
Operator splitting algorithms are a cornerstone of modern first-order optimization, relying critically on proximal operators as their fundamental building blocks. However, explicit formulas for proximal operators are available only for limited classes of functions, restricting the applicability of these methods. Recent work introduced HJ-Prox \cite{OsherHeatonFung2023}, a zeroth-order Monte Carlo approximation of the proximal operator derived from Hamilton–Jacobi PDEs, which circumvents the need for closed-form solutions. 
In this work, we extend the scope of HJ-Prox by establishing that it can be seamlessly incorporated into operator splitting schemes while preserving convergence guarantees. In particular, we show that replacing exact proximal steps with HJ-Prox approximations in algorithms such as proximal gradient descent, Douglas–Rachford splitting, Davis–Yin splitting, and the primal–dual hybrid gradient method still ensures convergence under mild conditions.
\end{abstract}

\begin{keywords}%
 Proximal, Operator Splitting, Derivative-Free, Zeroth-Order, Optimization, Monte Carlo, Hamilton-Jacobi
\end{keywords}

\section{Introduction}
In modern machine learning and optimization, splitting algorithms play an important role in solving complex problems, particularly those with nonsmooth composite objective functions \cite{ParikhBoyd2014}. 
Splitting algorithms face difficulty when a step involving the proximal operator lacks a closed-form solution, calling for computationally expensive and complex inner iterations to solve sub-optimization problems \cite{Tibshirani2017, TibshiraniTaylor2011}. To address this challenge, we build on recent work of Osher, Heaton, and Wu Fung, who showed that Hamilton–Jacobi (HJ) equations can be used to approximate proximal operators via a Monte Carlo scheme, termed HJ-Prox.~\cite{OsherHeatonFung2023}. In this work, we propose a new framework for splitting algorithms that replace the exact proximal operator with the HJ-Prox approximation. Our primary contribution is a theoretical and empirical demonstration that this new general framework maintains convergence near the true solution, reducing the need for proximal calculus and introducing a more universal and readily applicable approach to splitting algorithms. For this workshop paper, we focus on proximal gradient descent (PGD), Douglas Rachford Splitting (DRS), Davis-Yin Splitting (DYS), and the primal-dual hybrid gradient algorithm (PDHG)~\cite{ryu2022large}. 

\section{Background}
Splitting algorithms are designed to solve composite convex optimization problems of the form
\begin{equation}
    \min_x f(x) + g(x),
    \label{eq:min_prob}
\end{equation}
$f$ and $g$ are proper, lower‑semicontinuous (LSC) and convex. Their efficiency, however, depends critically on the availability of closed-form proximal operators. When these operators are unavailable, the proximal step must be approximated through iterative subroutines, creating a substantial computational bottleneck.
To address this challenge, several lines of research have emerged. One approach focuses on improving efficiency through randomization within the algorithmic structure. These methods reduce computational cost by sampling blocks of variables, probabilistically skipping the proximal step, or solving suboptimization problems incompletely with controlled error \cite{Mishchenko2022ProxSkip, BonettiniPratoRebegoldi2020OO, BricenoAriasEtAl2019, CondatRichtarik2022RandProxOO}. Alternatively, other approaches reformulate the problem by focusing on dual formulations \cite{Tibshirani2017, MazumderHastie2012}.

While these techniques improve efficiency, they share common limitations. They require extensive derivations and complex analysis to handle the proximal operator, and proximal operators remain problem-dependent, typically requiring tailored solution strategies for each specific function class. This creates a critical research gap: the need for a generalizable method that can approximate the proximal operator without derivative information, making it suitable for zeroth-order optimization problems where only function evaluations are available.

\subsection{Hamilton-Jacobi-based Proximal (HJ-Prox)}

A promising solution to this challenge has emerged from recent work that approximates the proximal operator using a Monte-Carlo approach inspired by Hamilton-Jacobi (HJ) PDEs. Specifically, Osher, Heaton, and Wu Fung \cite{OsherHeatonFung2023} showed that 
\begin{align}
    \prox_{tf}(x) & \amp = \amp \lim_{\delta \to 0^+}\frac{\mathbb{E}_{y \sim \mathnormal{\mathcal{N}(x,\delta t I)}}\left[y \cdot\exp(-f(y)/\delta)\right]}{\mathbb{E}_{y \sim \mathnormal{\mathcal{N}(x,\delta t I)}}\left[\exp(-f(y)/\delta)\right]}
    \\
    &\amp \approx \amp \frac{\mathbb{E}_{y \sim \mathnormal{\mathcal{N}(x,\delta t I)}}\left[y \cdot\exp(-f(y)/\delta)\right]}{\mathbb{E}_{y \sim \mathnormal{\mathcal{N}(x,\delta t I)}}\left[\exp(-f(y)/\delta)\right]} \quad \text{ for some $\delta > 0$}
    \\
    & \amp = \amp \hjprox_{tf}(x)
    \label{eq:hj_prox}
\end{align}
where $\mathcal{N}(x, \delta t I)$ represents the normal distribution with mean $x$ and covariance matrix $\delta t I$, $t > 0$, and $f$ is assumed to be weakly-convex~\cite{ryu2022large}. 

The HJ-Prox, denoted by $\hjprox_{tf}$ in~\eqref{eq:hj_prox}, fixes a small value of $\delta > 0$ to approximate the limiting expression above, enabling a Monte Carlo approximation of the proximal operator in a zeroth-order manner~\cite{OsherHeatonFung2023, tibshirani2025laplace, heaton2024global, meng2025recent}. 

This approach is particularly attractive because it requires only function evaluations, avoiding the need for derivatives or closed-form solutions.
Subsequent research has explored HJ-Prox applications, primarily in global optimization via adaptive proximal point algorithms \cite{heaton2024global, ZhangEtAl2024, zhang2025thinking}. However, these applications have remained narrow in scope, focusing on specific algorithmic contexts rather than establishing a general framework.
Our work expands upon the theory of HJ-Prox by creating a comprehensive framework that can be applied to the entire family of splitting algorithms for convex optimization, including proximal gradient descent (PGD)~\cite{rockafellar1970convex, ryu2022large}, Douglas Rachford Splitting (DRS)~\cite{lions1979splitting, eckstein1992douglas}, Davis-Yin Splitting (DYS)~\cite{DavisYin2015}, and primal-dual hybrid gradient (PDHG)~\cite{ChambollePock2011}. To our knowledge, the direct approximation of the proximal operator via HJ equations for use in general splitting methods has not been previously explored.

\section{HJ-Prox-based Operator Splitting}
We now show how HJ-Prox can be incorporated into splitting algorithms such as PGD, DRS, DYS, and PDHG. The key idea is simple: by replacing exact proximal steps with their HJ-Prox approximations, we retain convergence guarantees while eliminating the need for closed-form proximal formulas or costly inner optimization loops. For readability, all proofs are deferred to the Appendix.

Our analysis builds on a classical result concerning perturbed fixed-point iterations. In particular, Combettes \cite[Thm. 5.2]{Combettes2001} established convergence of Krasnosel’skiĭ–Mann (KM) iterations subject to summable errors:
\newcommand{\perturbedKMTheorem}[1]{
    Let $\left\{x_k\right\}_{k\geq0}$ be an arbitrary sequence generated by
    \begin{equation}
        x_{k+1} \amp = \amp x_k + \lambda_k\left({T}x_k + \epsilon_k - x_k\right), 
    \end{equation}
    where $T \colon \Real^n \to \Real^n$ is an operator that has at least one fixed point.
    If $\{\|\epsilon_k\|\}_{k\geq0} \in \ell^1$ (that is, $\epsilon_k$ is summable), $T - I$ is demiclosed at $0$, and  $\{\lambda_k\}_k \geq 0$ lies in $[\gamma, 2- \gamma]$ for some $\gamma \in (0,1)$, then $\{x_k\}_{k\geq0}$ converges to a fixed point of $T$. 
    }
 \begin{theorem}[Convergence of Perturbed Krasnosel'skii-Mann Iterates]
    \label{theorem:perturbed_KM}
    \perturbedKMTheorem{main}
\end{theorem}
Thus, to establish convergence of HJ-Prox–based splitting, it suffices to bound the HJ approximation error. The following result, originally proved in \cite{OsherHeatonFung2023, ZhangEtAl2024}, provides the required bound.

\newcommand{\hjproxerrorTheorem}[1]{
Let $f : \Real^n \mapsto \Real$ be LSC. Then the Hamilton-Jacobi approximation incurs errors that are uniformly bounded.
    \begin{eqnarray}
    \sup_{x} \left\|\hjprox_{tf}(x) - \prox_{tf}(x)\right\| & = & \sqrt{2n t \delta}.
    \end{eqnarray}
}
\begin{theorem}[Error Bound on HJ-Prox]
\label{thm:hj_prox_error_bound}
\hjproxerrorTheorem{main}
\end{theorem}
This uniform error bound guides the choice of $\delta$ in each iteration of our splitting algorithms. In particular, by selecting $\delta_k$ so that the resulting error sequence is summable, Theorem \ref{theorem:perturbed_KM} ensures convergence of the HJ-Prox–based methods.

\newcommand{\pgdTheorem}[1]{ 
Let $f,g$ be proper, LSC, and convex, with $f$ additionally $L$-smooth.
    Consider the HJ-Prox-based PGD iteration given by 
    \begin{eqnarray}
        x_{k+1} & = & \prox^{\delta_k}_{tg}\left(x_k - t \nabla f(x_k)\right), \quad k=1,\ldots, 
    \end{eqnarray}
    with step size $0 < t < 2/L$ and $\left\{ \sqrt{\delta_k} \right\}_{k\geq1}$ a  summable sequence.
    Then $x_k$ converges to a minimizer of $f + g$.
    }

We rely on Theorem~\ref{theorem:perturbed_KM} to prove the convergence of the four fixed-point methods that use HJ-Prox. 
For simplicity, take $\lambda_k = 1$ for all $k$. Recall that $T - I$ is demiclosed at 0 if $T$ is nonexpansive. In addition, recall that $T$ is nonexpansive if $T$ is averaged. Associated with each algorithm of interest is an algorithm map $T$ that takes the current iterate $x_k$ to the next iterate $x_{k+1}$. Consequently, to invoke Theorem~\ref{theorem:perturbed_KM}, we verify that $T$ is averaged and check the summability of the error introduced into $T$ by the HJ-prox approximation.
    
\begin{theorem}[HJ-Prox PGD]
    \label{theorem:hjpgd}
    \pgdTheorem{main}
\end{theorem}

\newcommand{\drsTheorem}[1]{
Let $f,g$ be proper, convex, and LSC. Consider the HJ-Prox–based DRS iteration given by  
    \begin{equation}
        \begin{split}
        x_{k+1/2} & \amp = \amp \prox^{\delta_k}_{tf}(z_k), \\
        x_{k+1} & \amp = \amp \prox^{\delta_k}_{tg}(2x_{k+1/2}-z_k), \\
        z_{k+1} & \amp = \amp z_k + x_{k+1} - x_{k+1/2},
        \end{split}
    \end{equation}
    with $\left\{ \sqrt{\delta_k} \right\}_{k\geq1}$ a summable sequence. Then $x_k$ converges to a minimizer of $f+g$.
}
\begin{theorem}[HJ-Prox DRS]
    \label{theorem:hjdrs}
    \drsTheorem{main}
\end{theorem}

\newcommand{\dysTheorem}[1]{
For DYS, consider $f + g + h$. Let $f,g,h$ be proper, LSC, and convex, with $h$ additionally $L$-smooth. Consider the HJ-Prox–based DYS algorithm given by
    \begin{equation}
        \begin{split}
        y_{k+1} & \amp = \amp \prox^{\delta_k}_{tf}(x_k), \\
        z_{k+1} & \amp = \amp \prox^{\delta_k}_{tg}\bigl(2y_{k+1} - x_k - t \nabla h(y_{k+1})\bigr),\\
        x_{k+1} & \amp = \amp x_k + z_{k+1} - y_{k+1},\\
        \end{split}
    \end{equation}
    with $\{ \sqrt{\delta_k}\}_{k\geq1}$ a summable sequence, and $0< t< 2/L$. Then $x_k$ converges to a minimizer of $f + g + h$.
}
\begin{theorem}[HJ-Prox DYS]
    \label{theorem:hjdys}
    \dysTheorem{main}
\end{theorem}

\newcommand{\pdhgTheorem}[1]{
Let $f,g$ be proper, convex, and LSC. Consider the HJ-Prox–based PDHG algorithm given by
    \begin{equation}
        \begin{split}
        y_{k+1} & \amp = \amp \prox^{\delta_k}_{\sigma g^*}(y_k + \sigma A x_k), \\
        x_{k+1} & \amp = \amp \prox^{\delta_k}_{\tau f}(x_k - \tau A^\top y_{k+1}),
        \end{split}
    \end{equation}
    with parameters $\tau,\sigma > 0$ satisfying $\tau\sigma \|A\|^2 < 1$ and $\{ \sqrt{\delta_k}\}_{k\geq1}$ a summable sequence. Where $g^*$ denotes the Fenchel conjugate of $g$. Then $x^k$ converges to a minimizer of $f(x)+g(Ax)$. 
}
\begin{theorem}[HJ-Prox PDHG]
    \label{theorem:hjpdhg}
    \pdhgTheorem{main}
\end{theorem}
\begin{figure}[H]\label{fig: LASSO and MLR}
  \centering
  \resizebox{\textwidth}{!}{%
  \begin{tabular}{@{}cccc@{}}
    \multicolumn{4}{c}{LASSO: $\underset{\beta}{\arg\min}\ \frac{1}{2}\|X\beta - y\|_2^2 + \lambda\|\beta\|_1$}\\
    \includegraphics[width=0.275\textwidth]{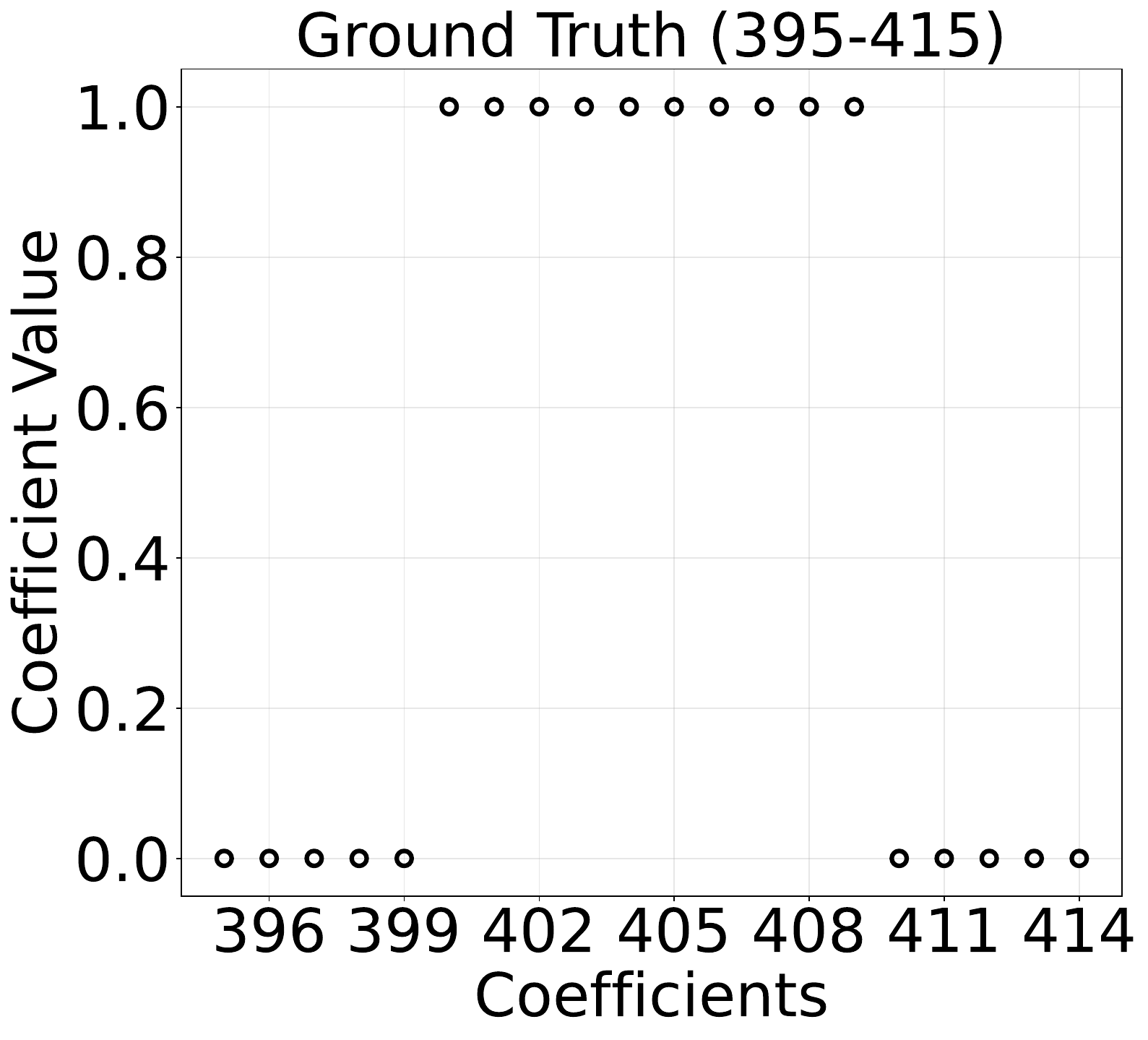} &
    \includegraphics[width=0.275\textwidth]{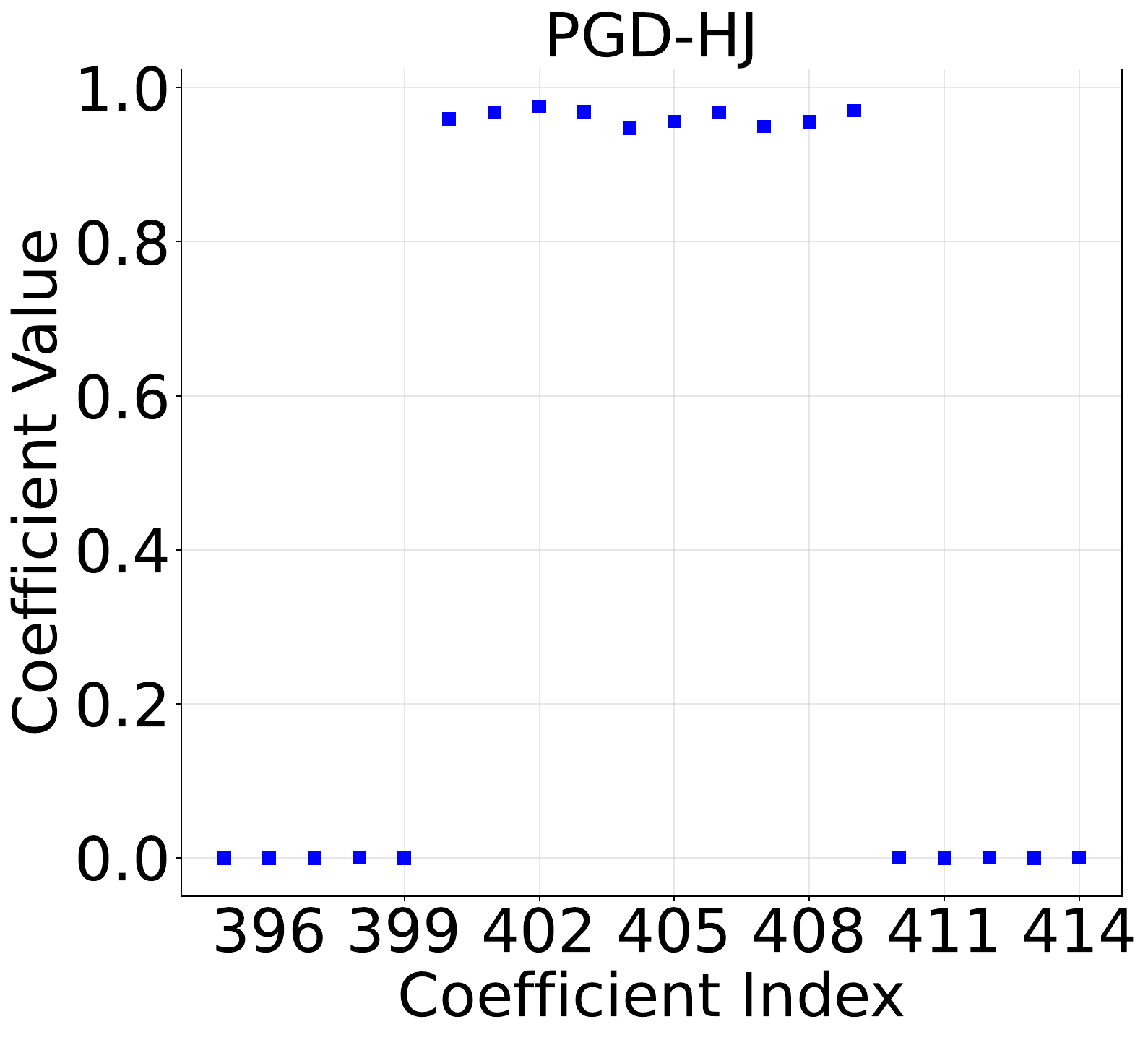} &
    \includegraphics[width=0.275\textwidth]{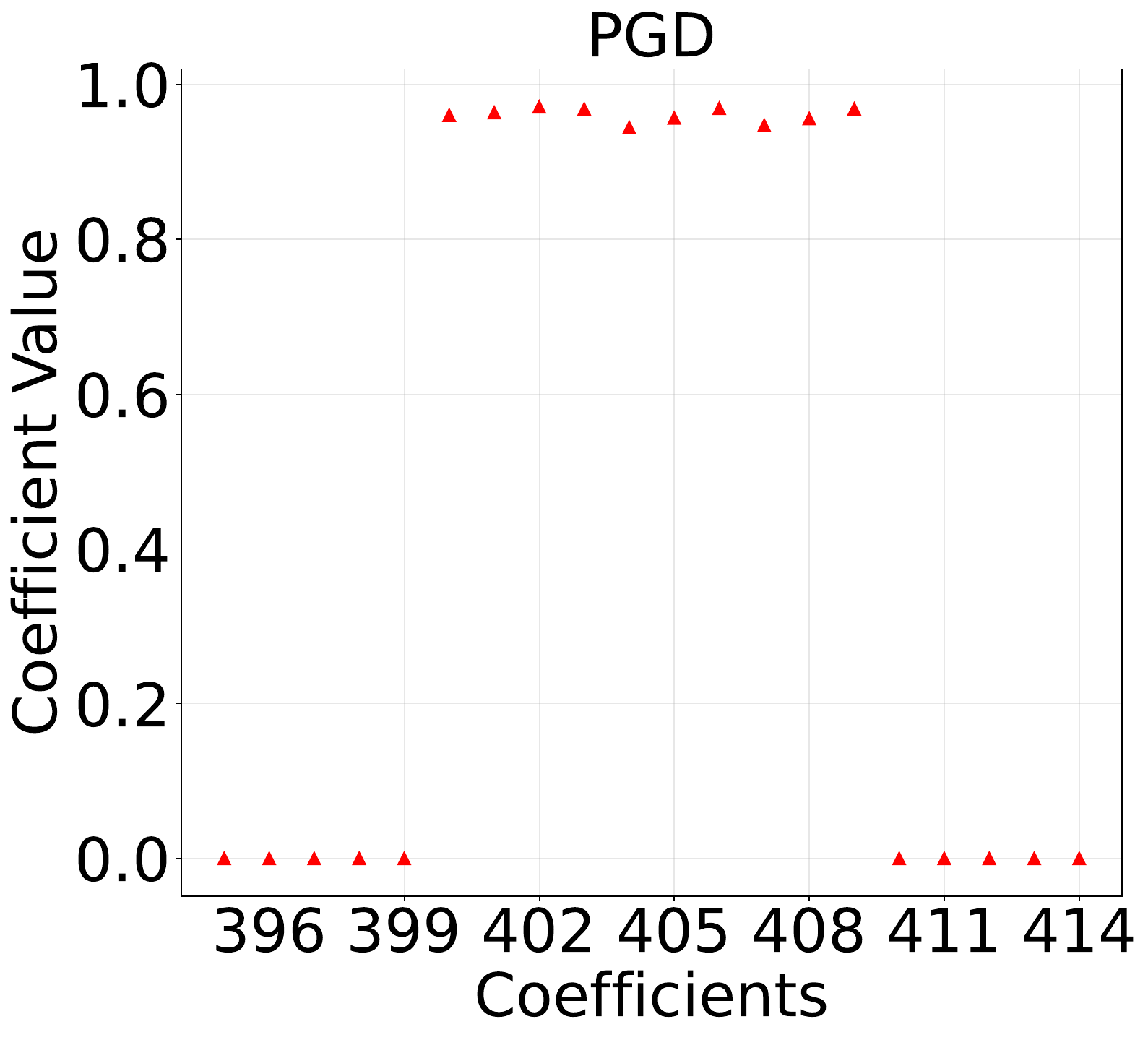} &
    \includegraphics[width=0.275\textwidth]{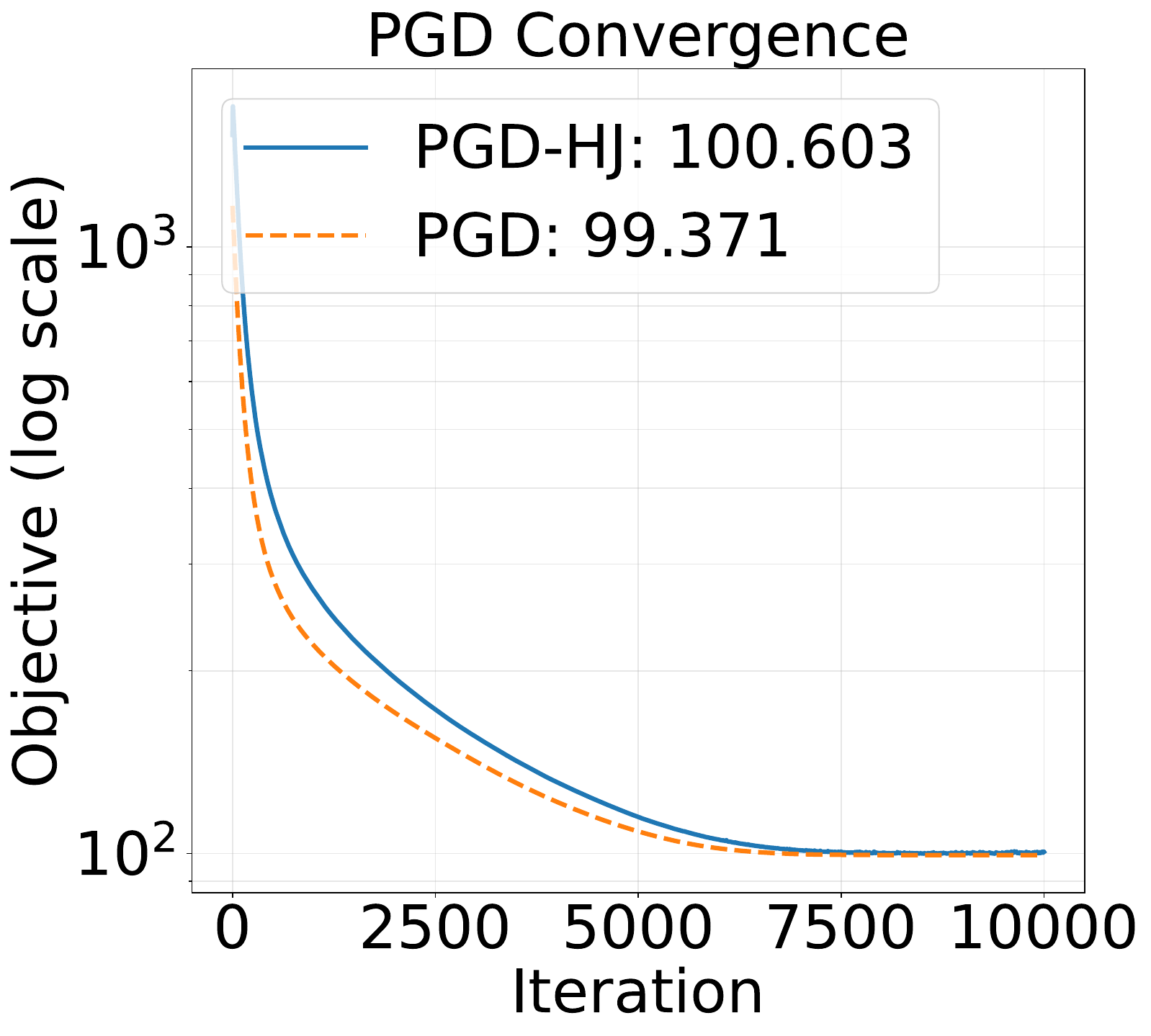} \\
    \multicolumn{4}{c}{Multitask Learning: $\underset{B}{\arg\min}\; \frac{1}{2}\,\| X B-Y\|_{F}^{2} +\lambda_{1}\,\| B\|_{*} +\lambda_{2}\sum_{i}\| b_{i,\cdot}\|_{2} +\lambda_{3}\sum_{j}\|b_{\cdot,j}\|_{2}$}\\
    \includegraphics[width=.275\textwidth]{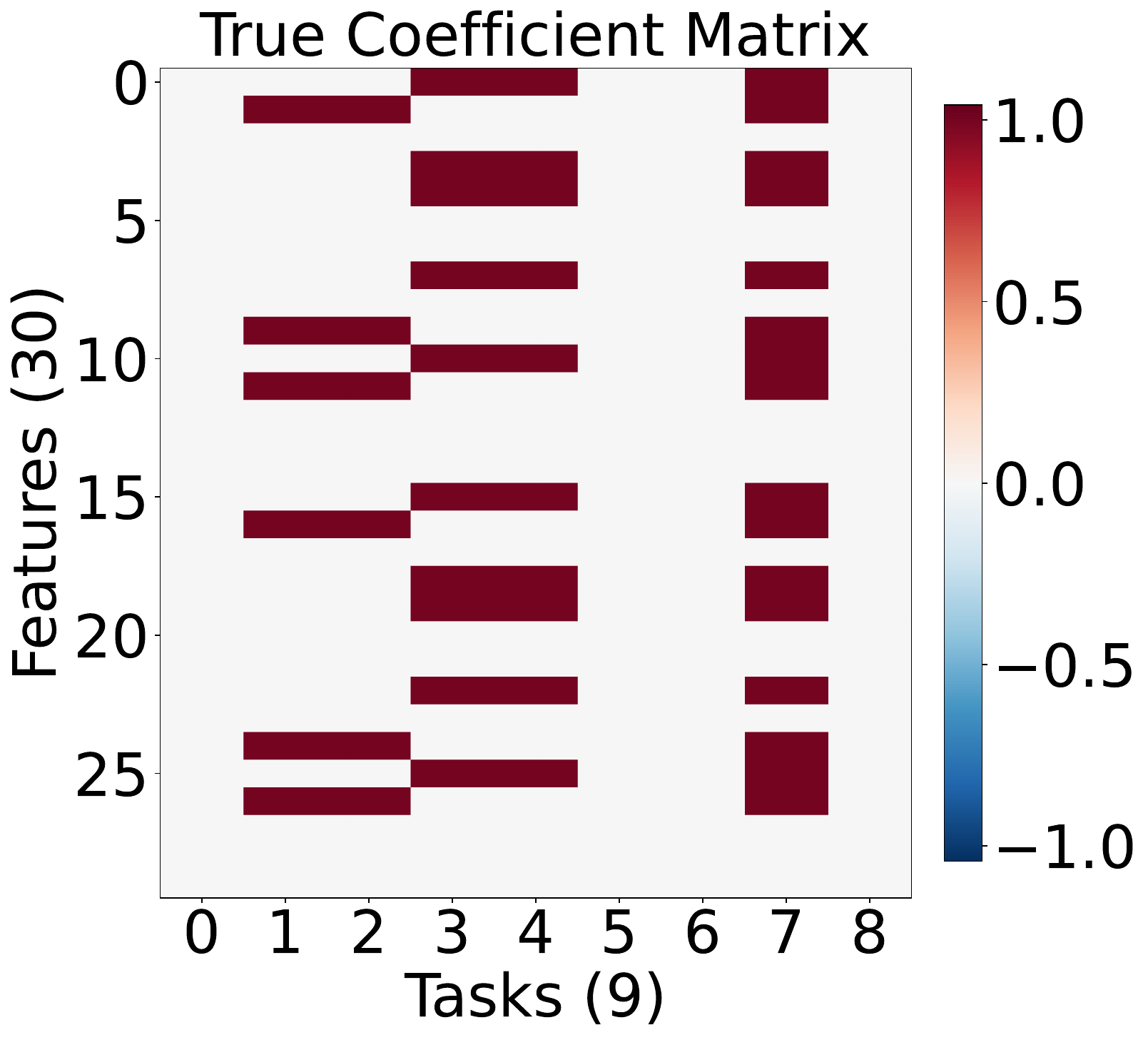} &
    \includegraphics[width=.275\textwidth]{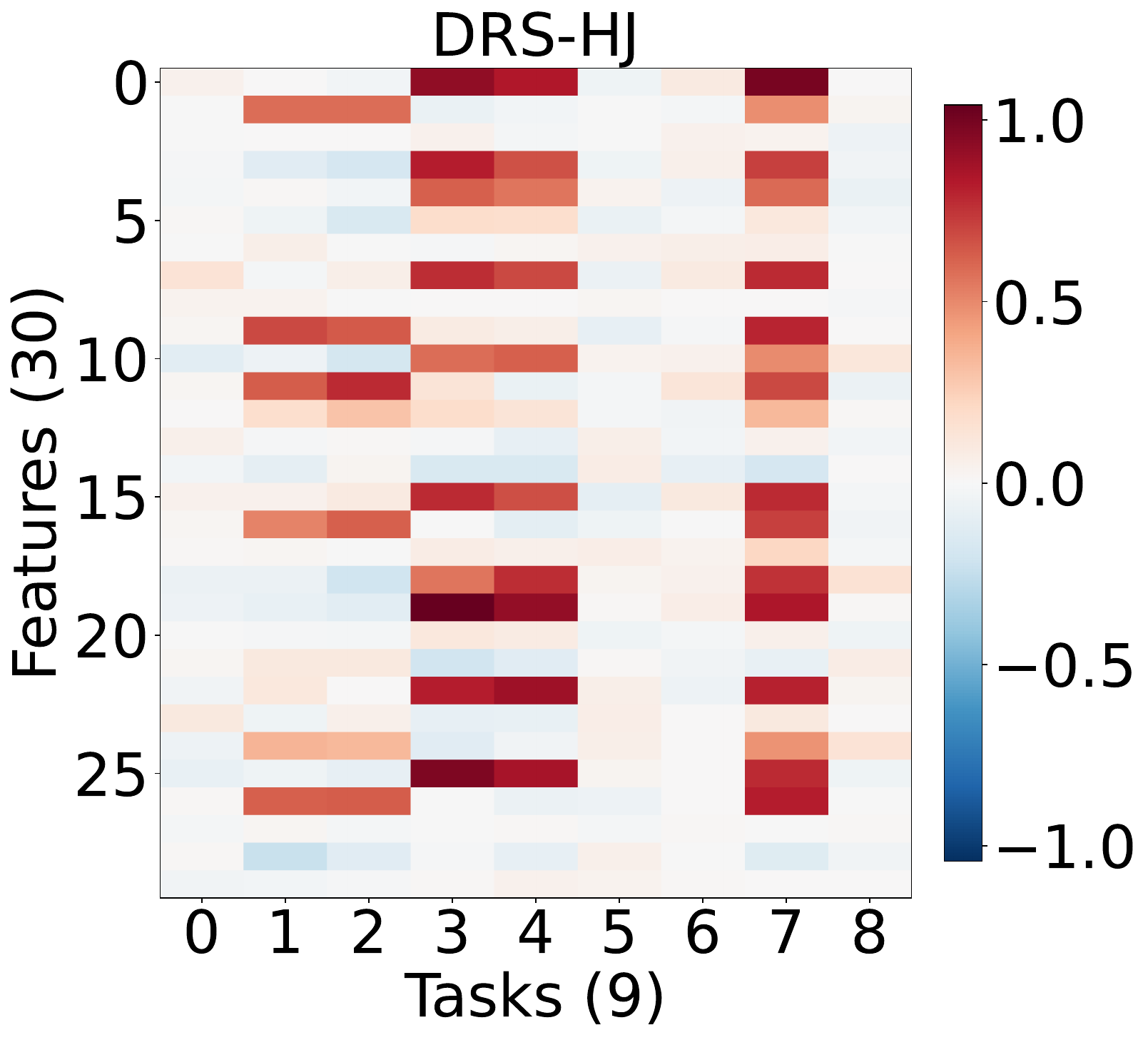} &
    \includegraphics[width=.275\textwidth]{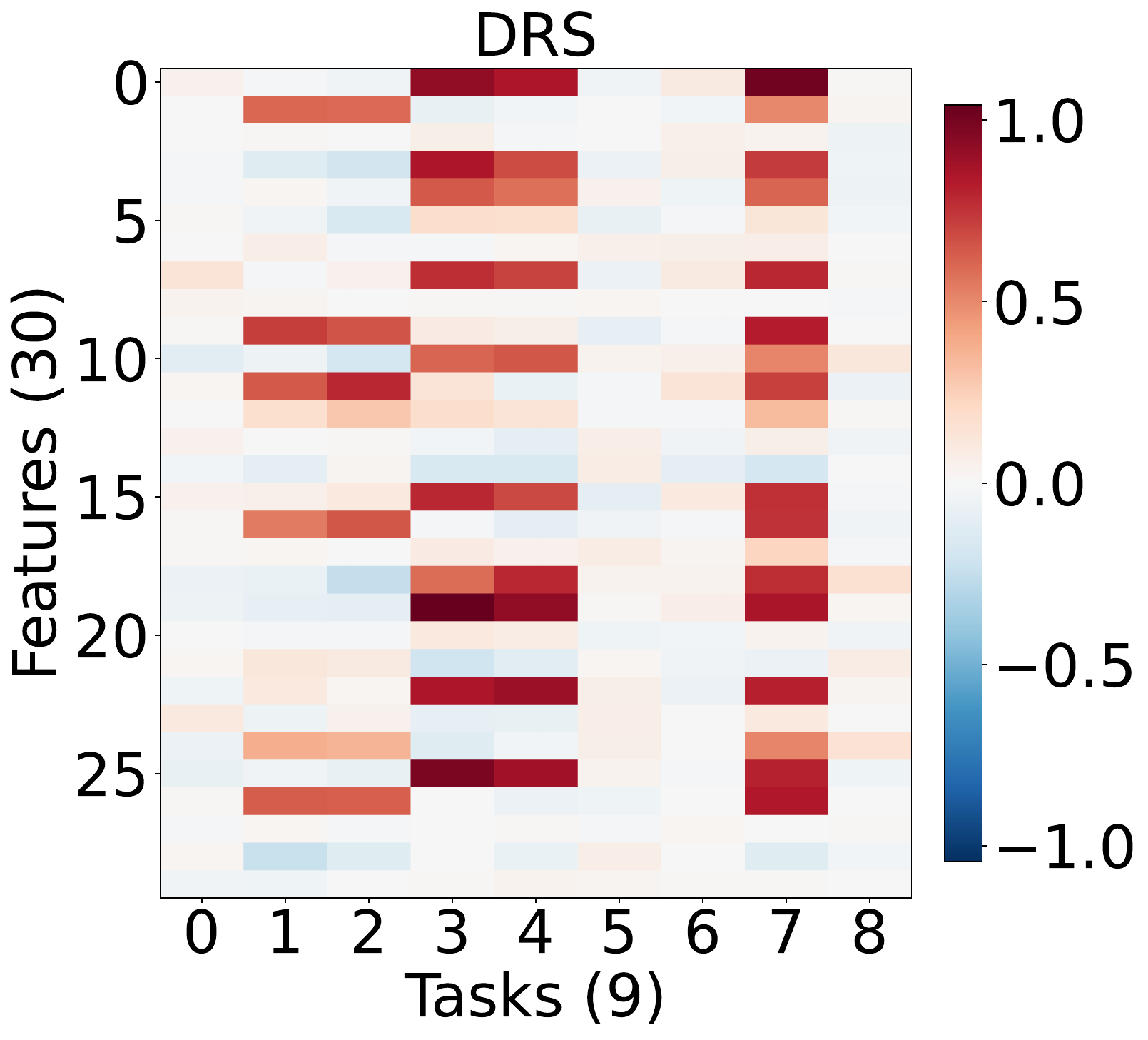} &
    \includegraphics[width=.275\textwidth]{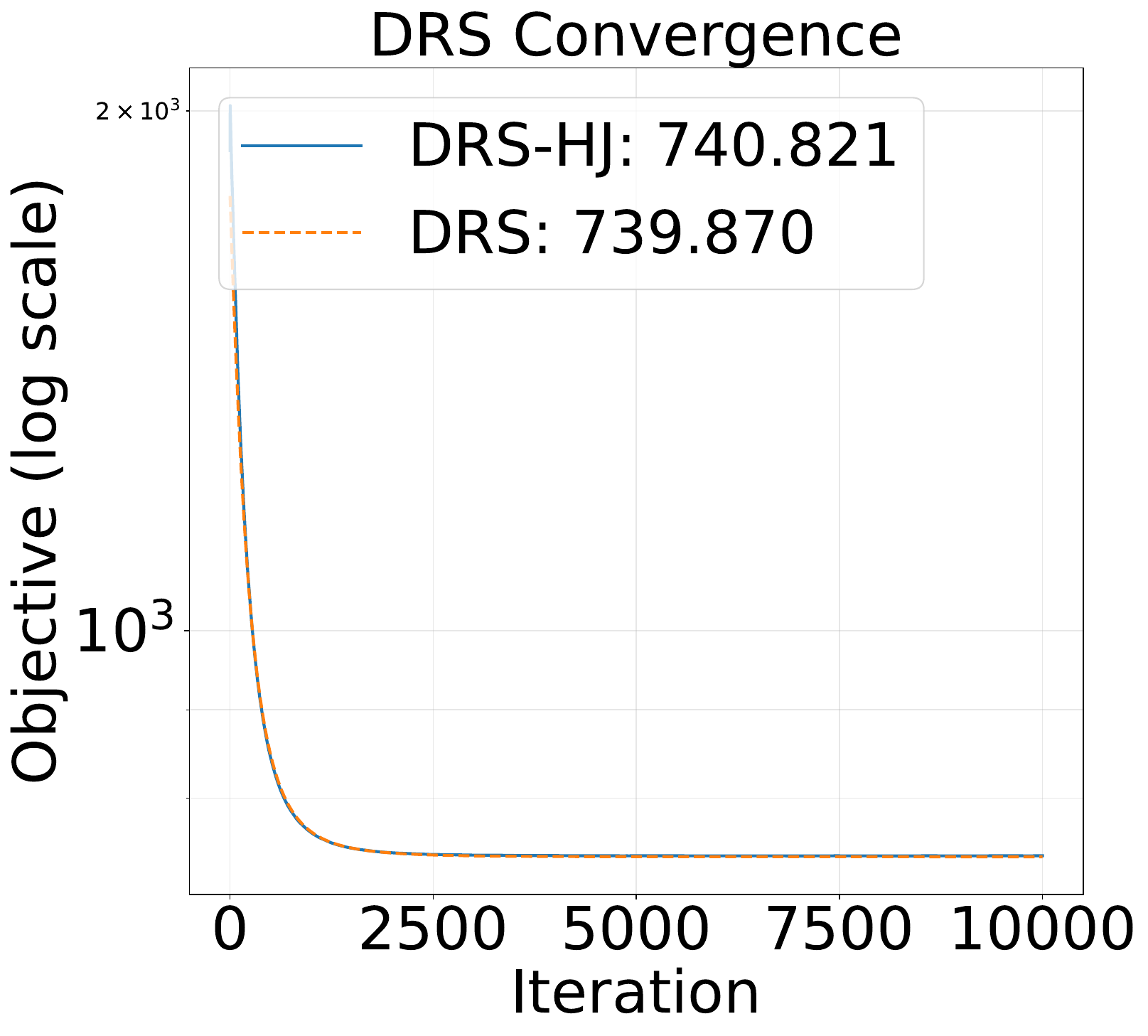}
  \end{tabular}%
  }
  \caption{LASSO and Multitask Learning Results}\label{fig: LASSO and MLR}
\end{figure}

\section{Experiments}
\begin{figure}[H]
    \centering
    \resizebox{\textwidth}{!}{%
    \begin{tabular}{@{}cccc@{}} 
        \multicolumn{4}{c}{Fused LASSO: $\underset{B}{\arg\min}\ \frac{1}{2}\|\beta - y\|^2 + \lambda\|D\beta\|_1$}
        \\
        \includegraphics[width=.275\textwidth]{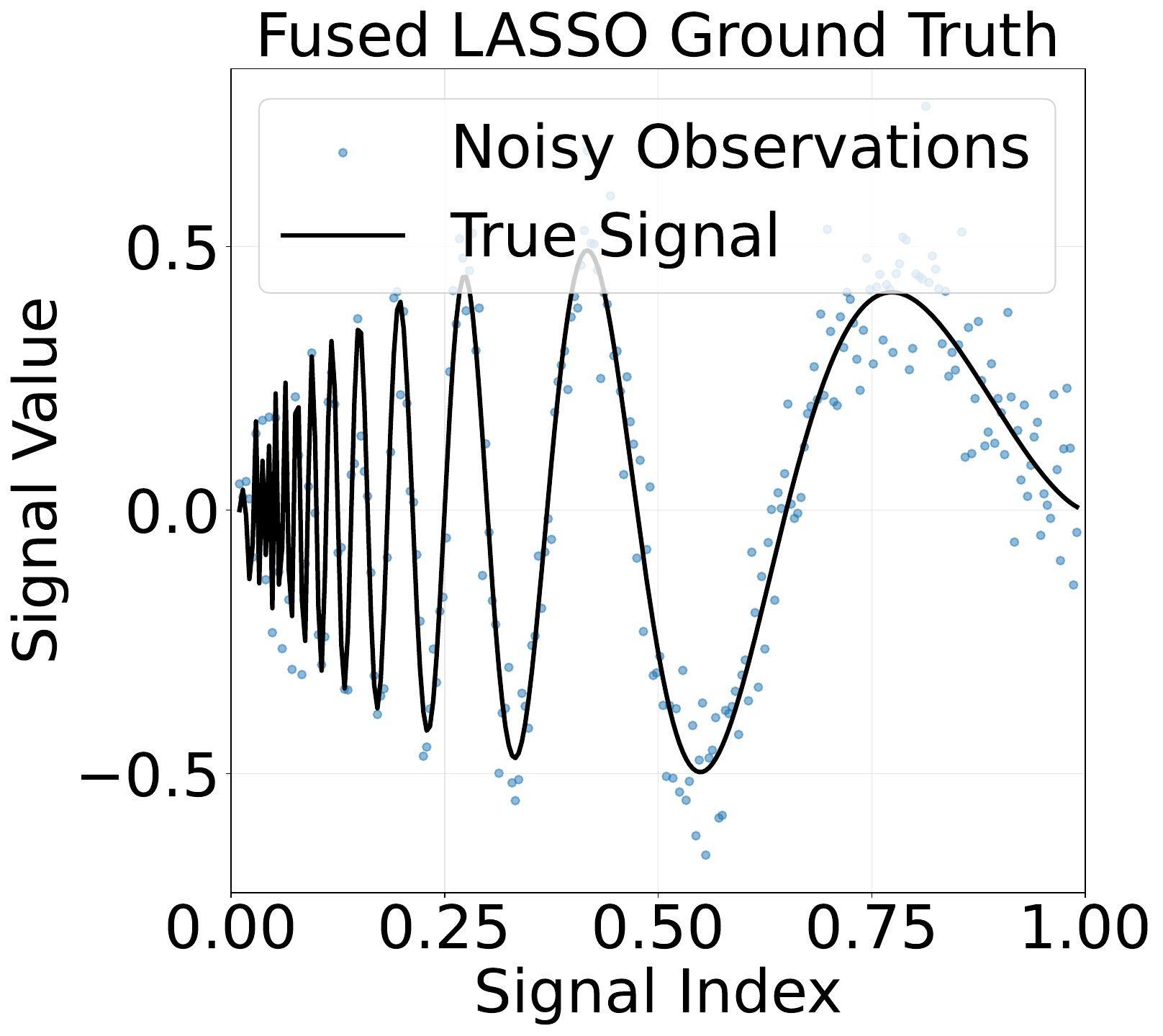}
        & 
        \includegraphics[width=.275\textwidth]{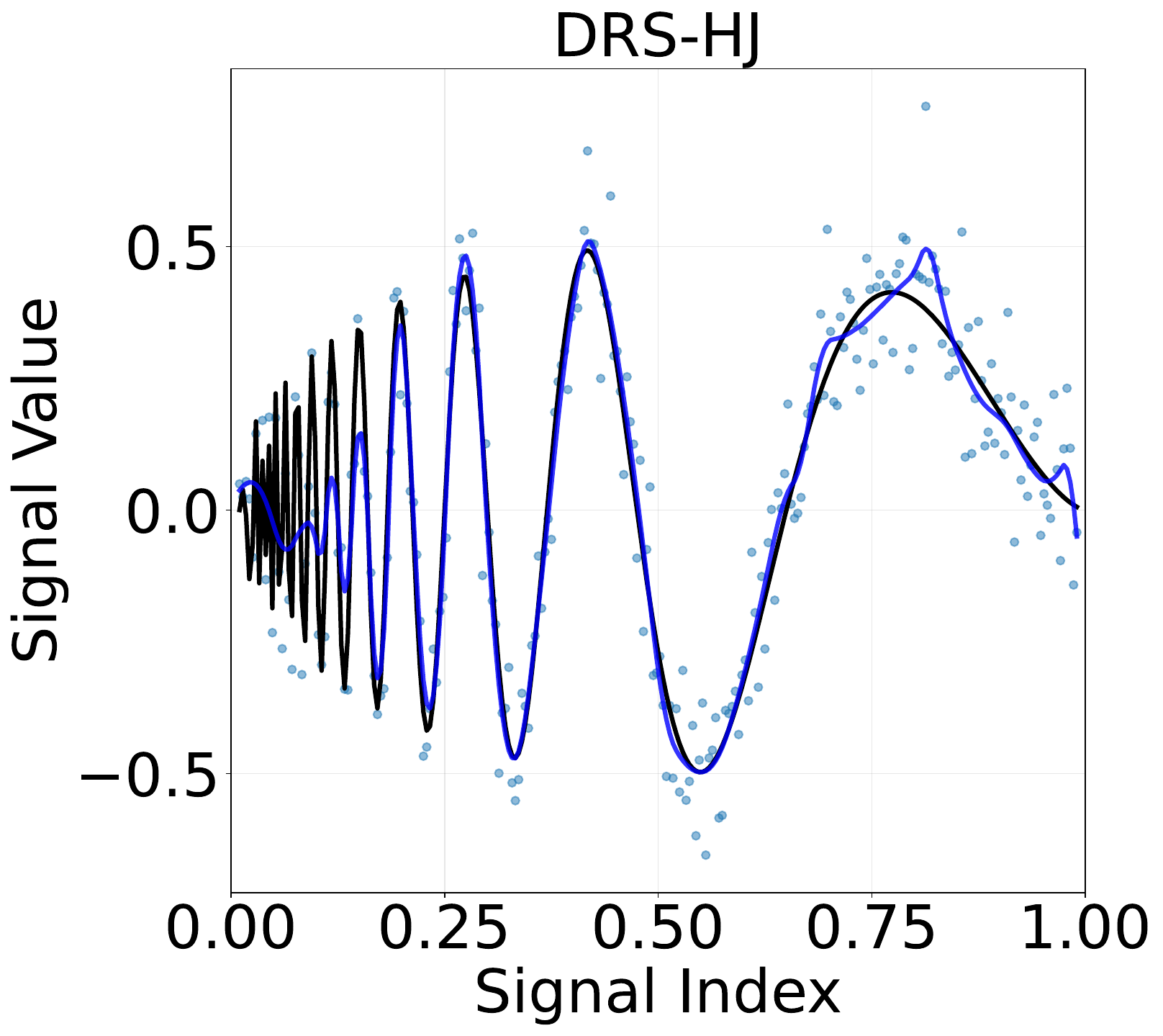}
        &
        \includegraphics[width=.275\textwidth]{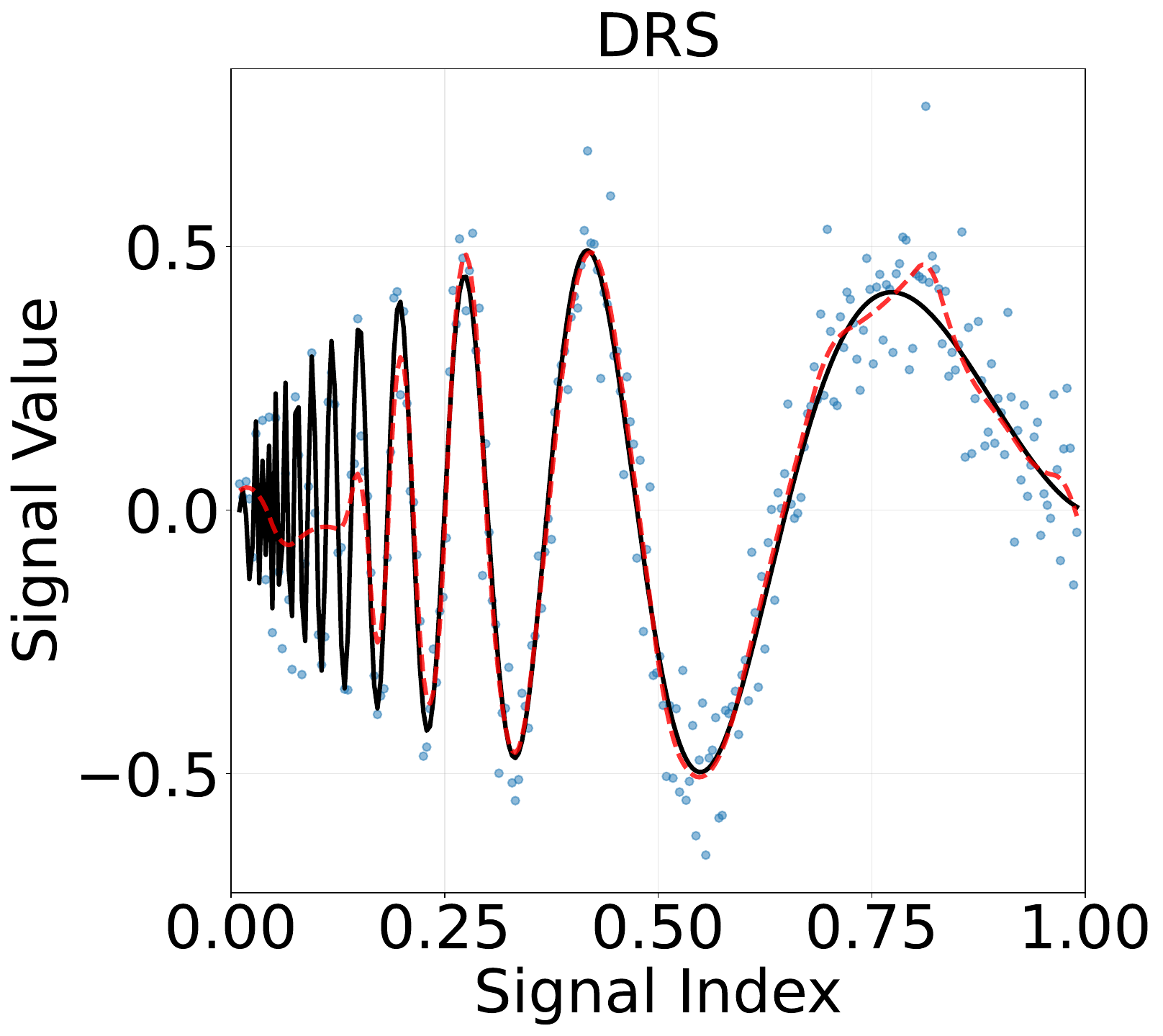}
        &
        \includegraphics[width=.275\textwidth]{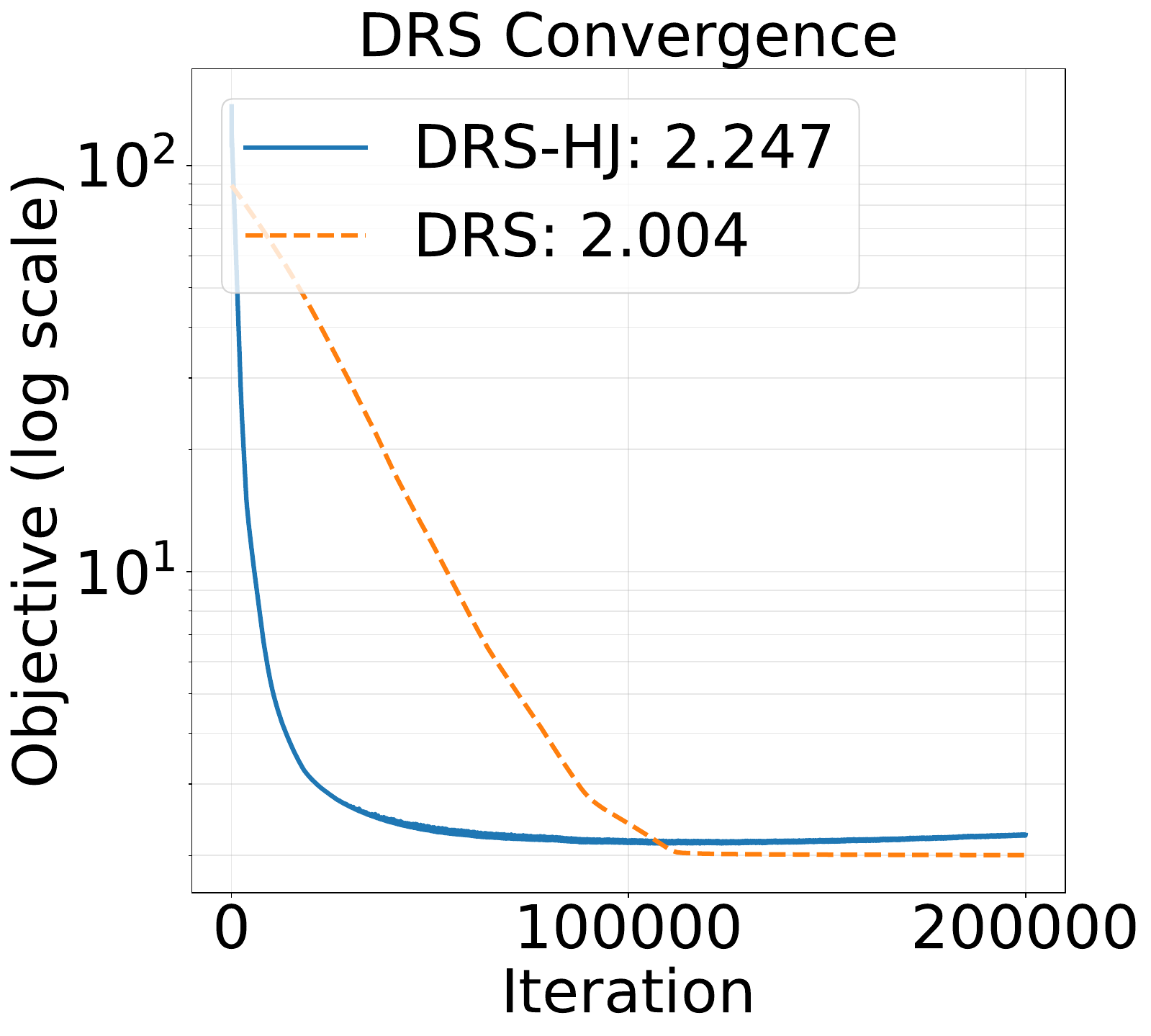} 
        \\
        \multicolumn{4}{c}{Sparse Group LASSO: $\underset{\beta}{\arg\min}\; \frac{1}{2}\,\| X \beta - y\|_{2}^{2} +\lambda_1\sum_{g=1}^{G}\| \beta_{g}\|_{2} +\lambda_2\| \beta\|_{1}$}
        \\
        \includegraphics[width=.275\textwidth]{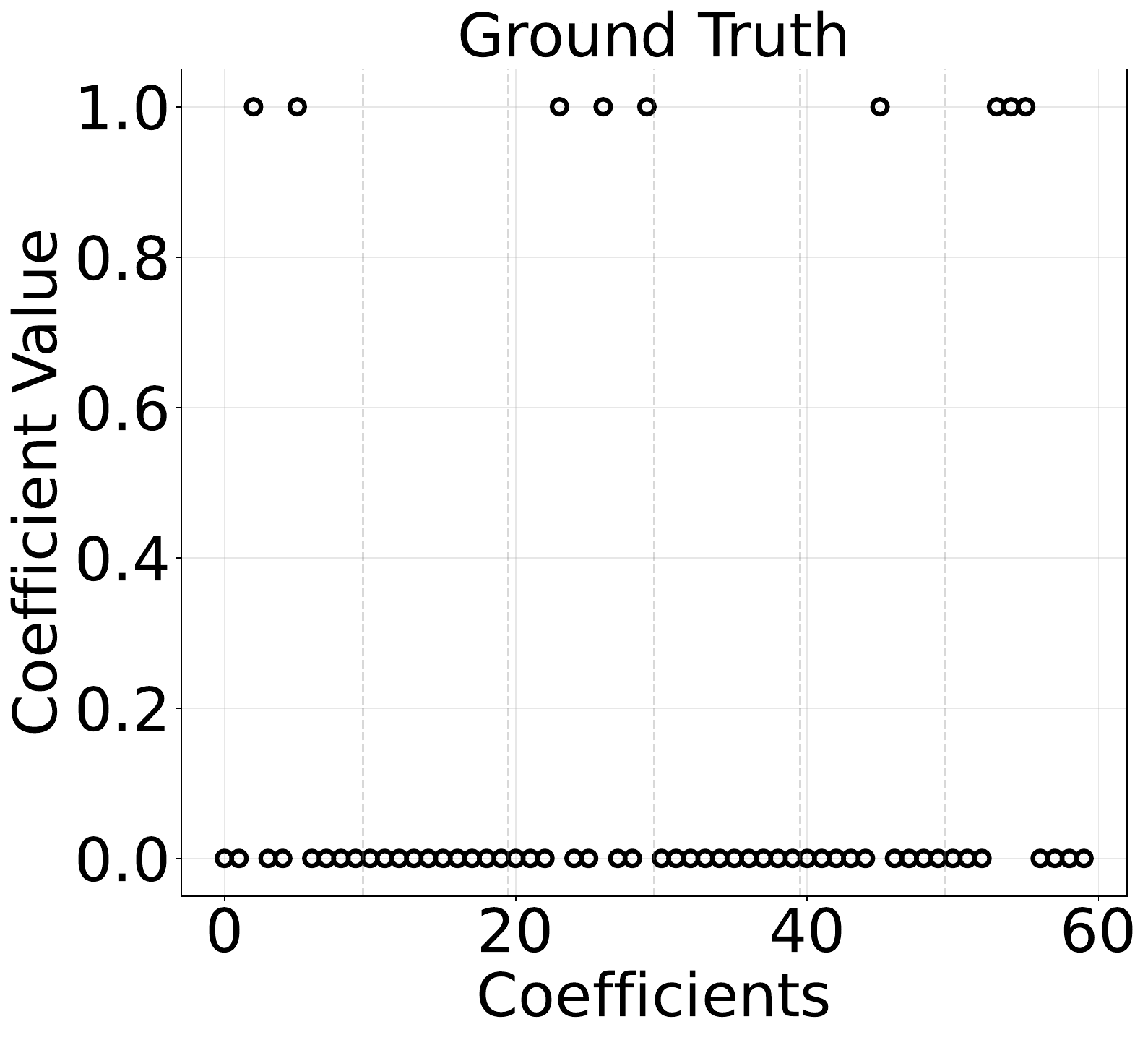}
        & 
        \includegraphics[width=.275\textwidth]{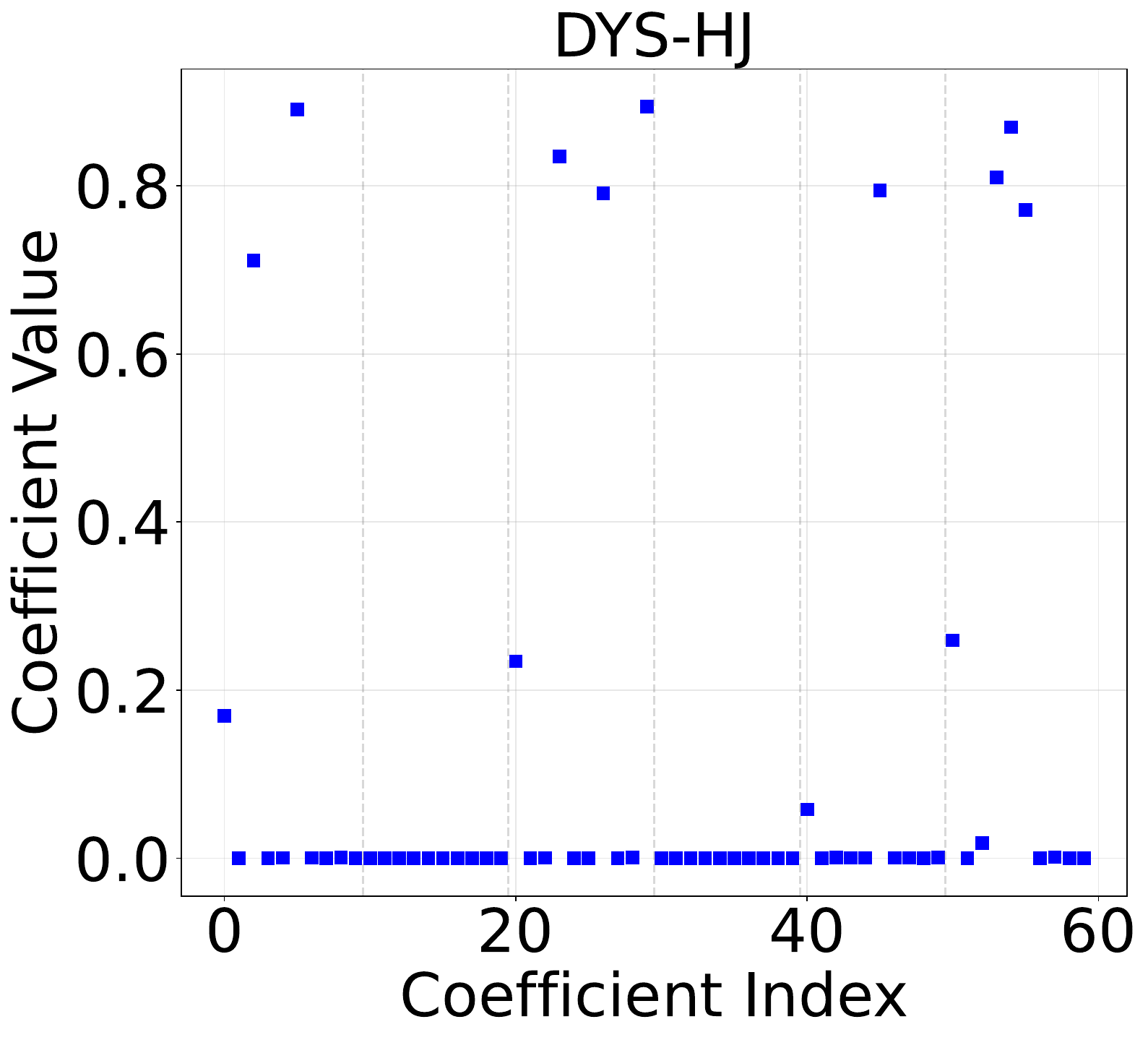}
        &
        \includegraphics[width=.275\textwidth]{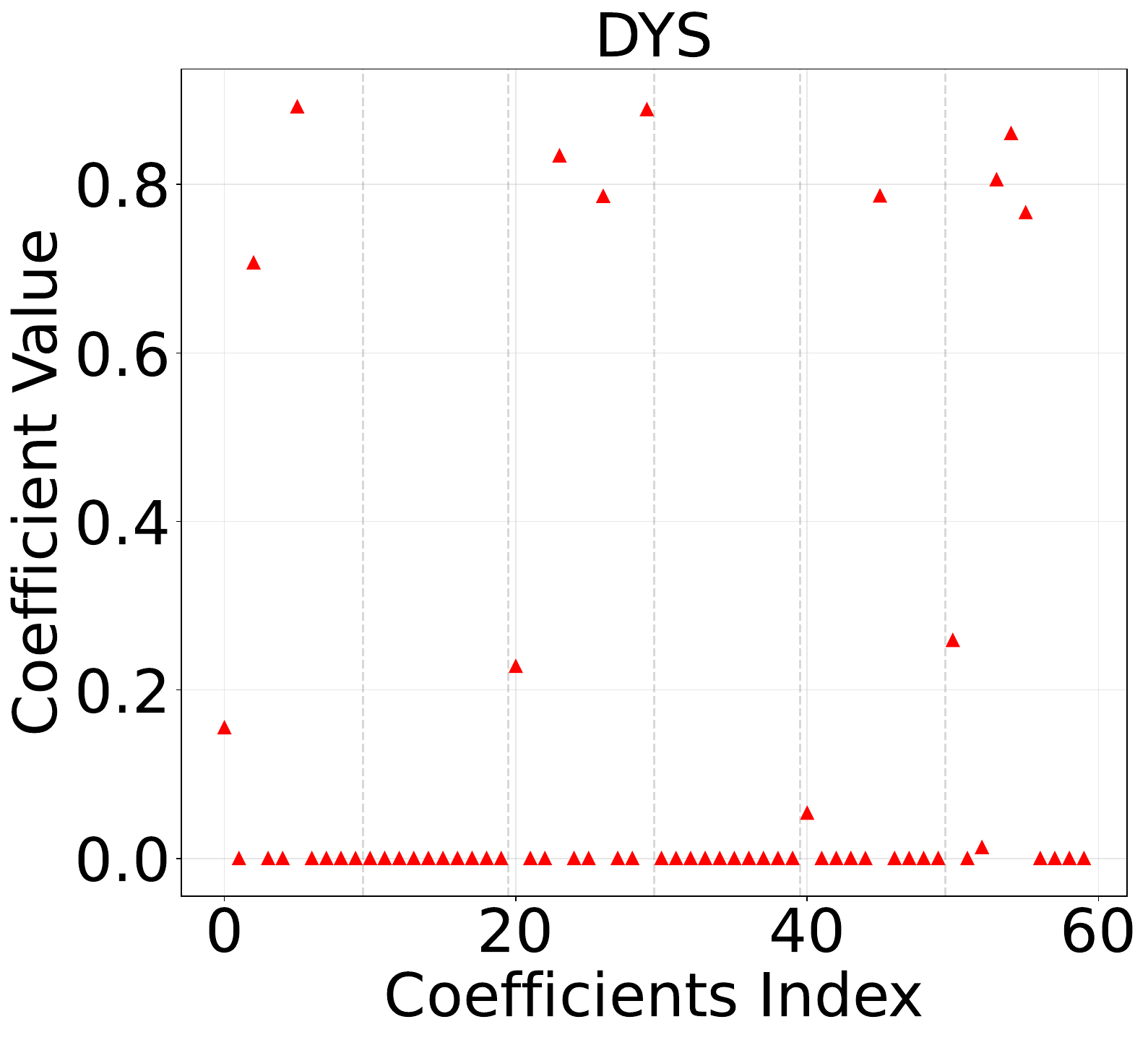}
        &
        \includegraphics[width=.275\textwidth]{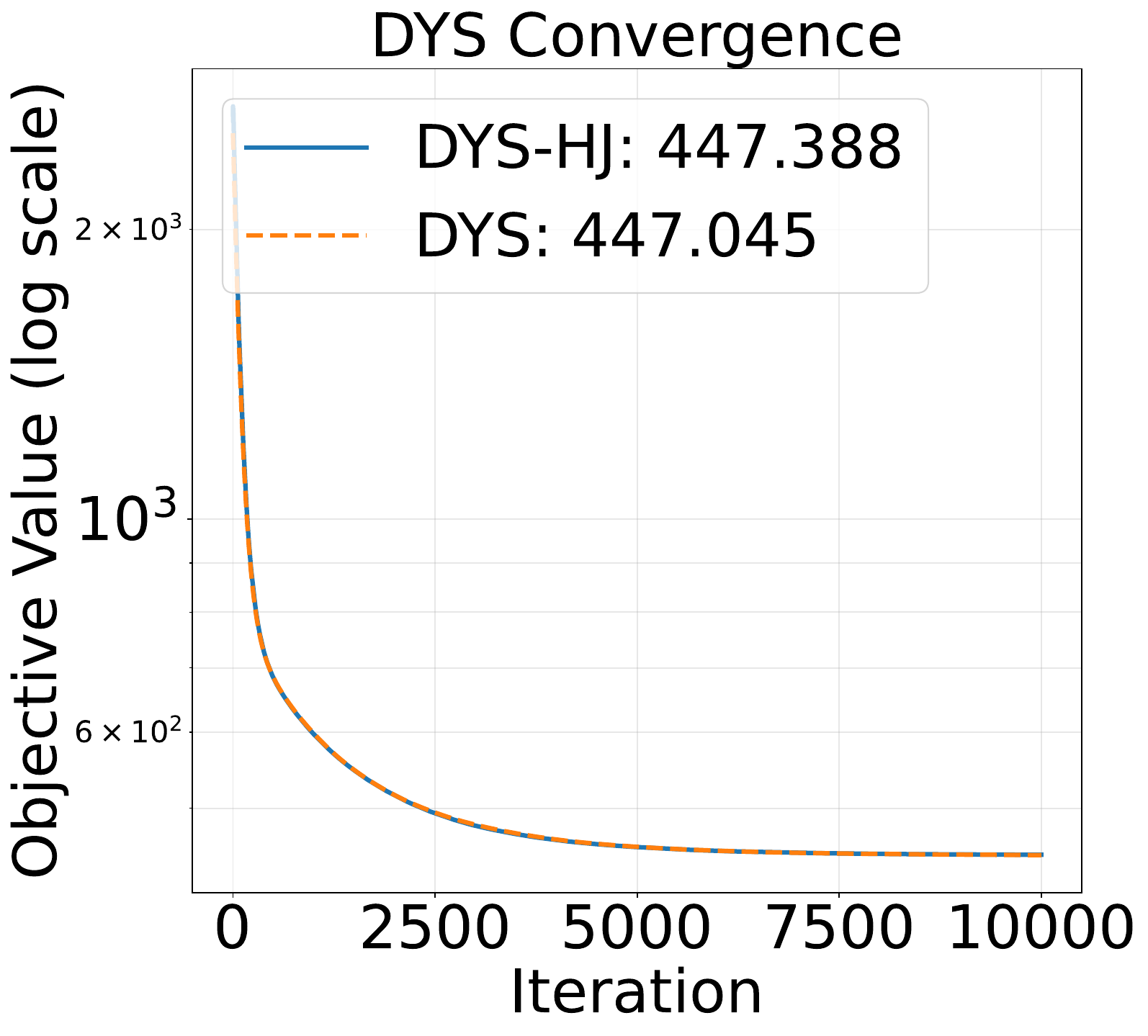}
    \end{tabular}%
    }
    \caption{Fused LASSO and Sparse Group LASSO Results}\label{fig: Fused LASSO and Sparse Group LASSO} 
\end{figure}

We conduct five experiments to assess the effectiveness of HJ-Prox integrated with proximal splitting methods. First, we use PGD to solve the LASSO, an $\ell_1$-regularized least-squares model for sparse feature selection displayed in figure \ref{fig: LASSO and MLR}. We then apply DRS to multitask learning for structured sparse matrix recovery and to a third-order fused-LASSO formulation on Doppler data that smooths the underlying signal by penalizing third-order finite differences. Results are shown in Figures \ref{fig: LASSO and MLR} and \ref{fig: Fused LASSO and Sparse Group LASSO}. Next, we solve the sparse group LASSO with DYS, which induces sparsity both across groups and within groups displayed in Figure \ref{fig: Fused LASSO and Sparse Group LASSO}. Finally, in Figure \ref{fig: pdhgResults}, we use the PDHG method for total-variation (TV) image denoising, preserving sharp edges by penalizing the $\ell_1$ norm of the image gradients on a noisy, blurred sample image. We use identical problem parameters for both the HJ-Prox and analytical methods. For each experiment, the HJ-Prox temperature parameter $\delta_k$ is scheduled to satisfy conditions for convergence established in Theorem \ref{theorem:perturbed_KM}. Each experiment is designed to visually compare respective recovered signals with the ground truth. We also report the convergence and last iteration of the objective function values in the legend for both approaches. Further experimental details are in the Appendix \ref{AP: Experiments}.
\subsection{Results}
\begin{figure}[H]
    \centering
    \resizebox{\textwidth}{!}{%
    \begin{tabular}{@{}cccc@{}} 
        \multicolumn{4}{c}{Total Variation Denoising: $\underset{\beta}{\arg\min}\; \frac{1}{2} \| X \beta - y\|_{F}^{2} +\lambda \text{TV}(\beta)$}
        \\
        \includegraphics[width=.275\textwidth]{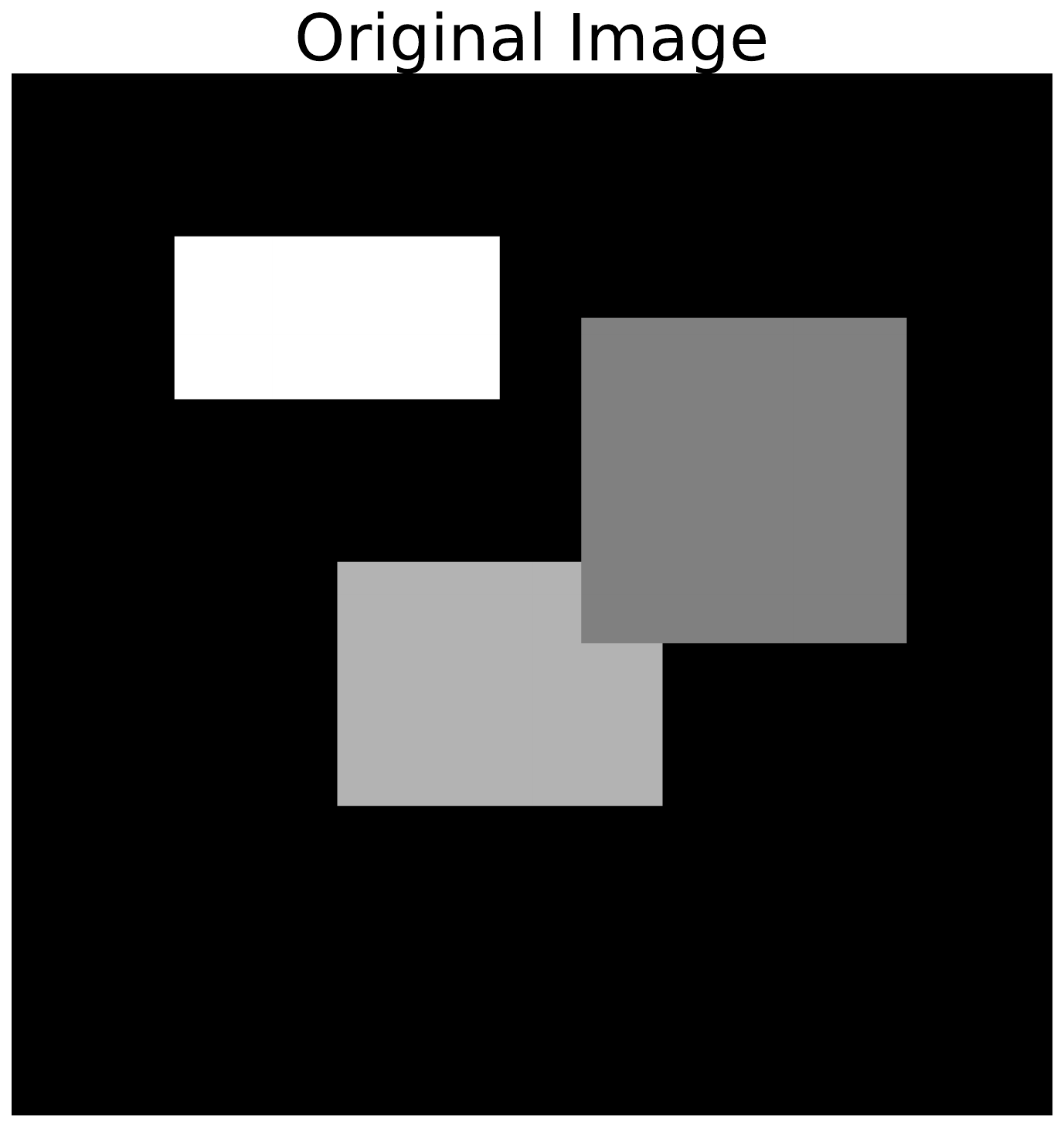}
        & 
        \includegraphics[width=.275\textwidth]{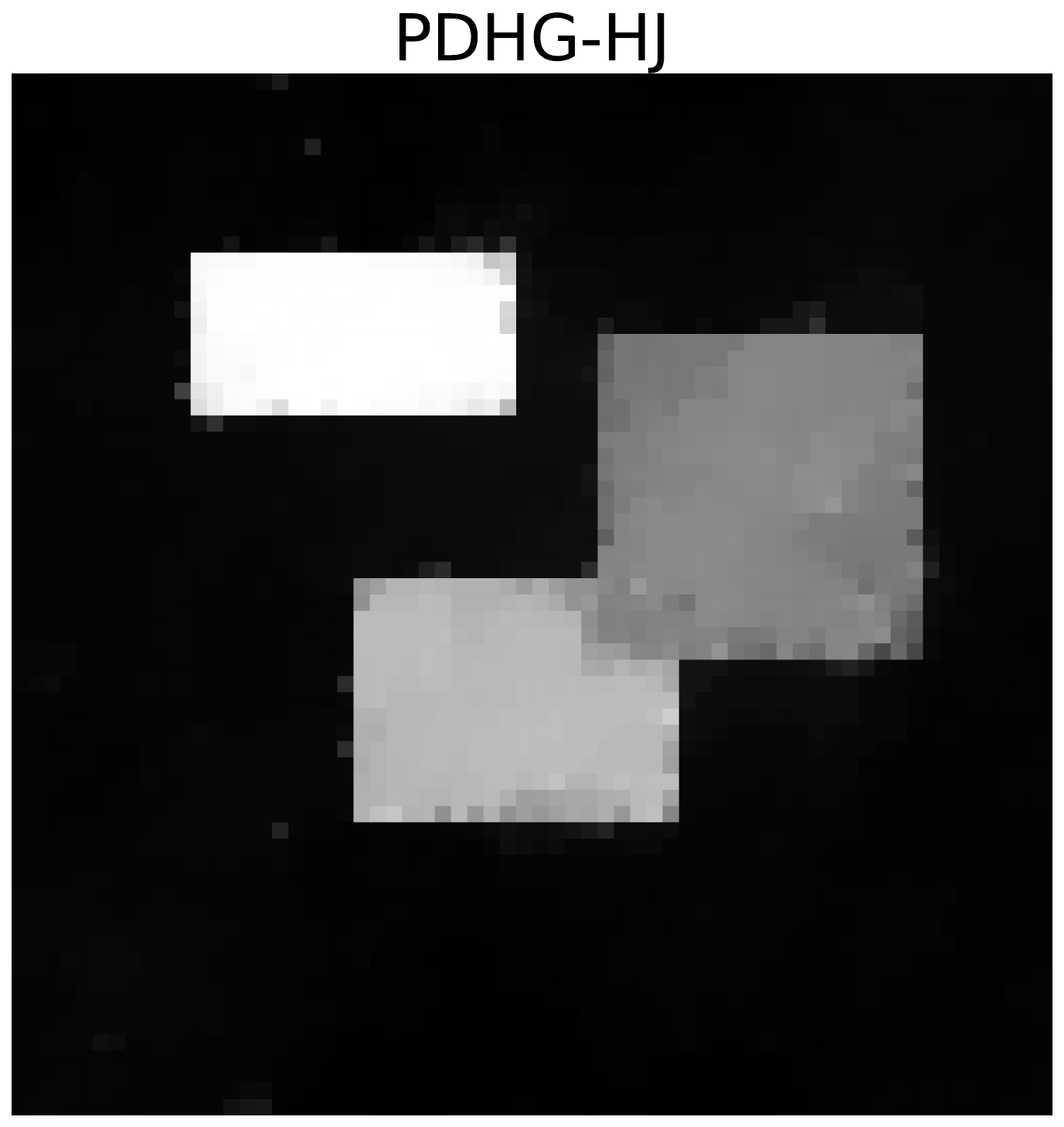}
        &
        \includegraphics[width=.275\textwidth]{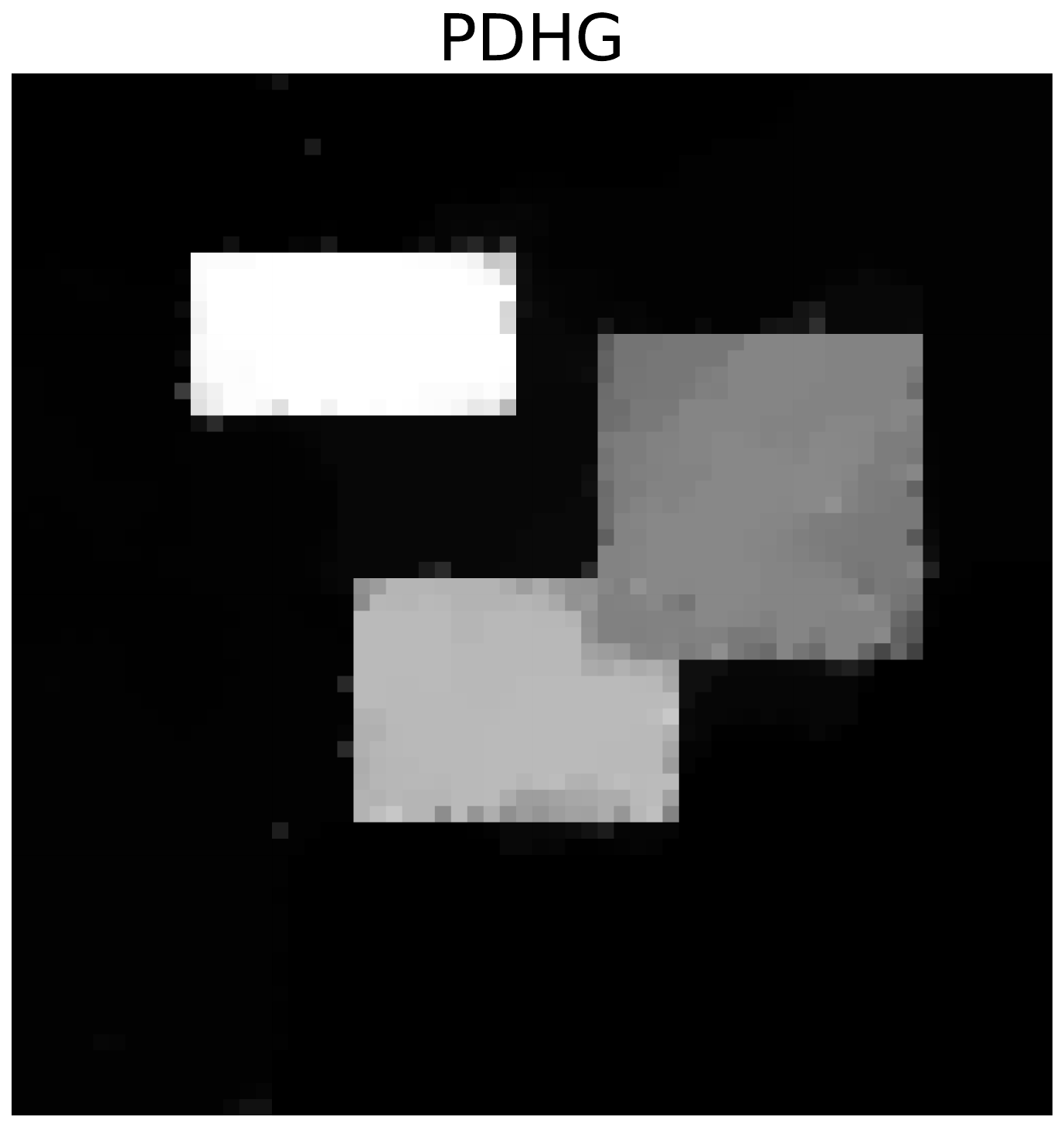}
        &
        \includegraphics[width=.275\textwidth]{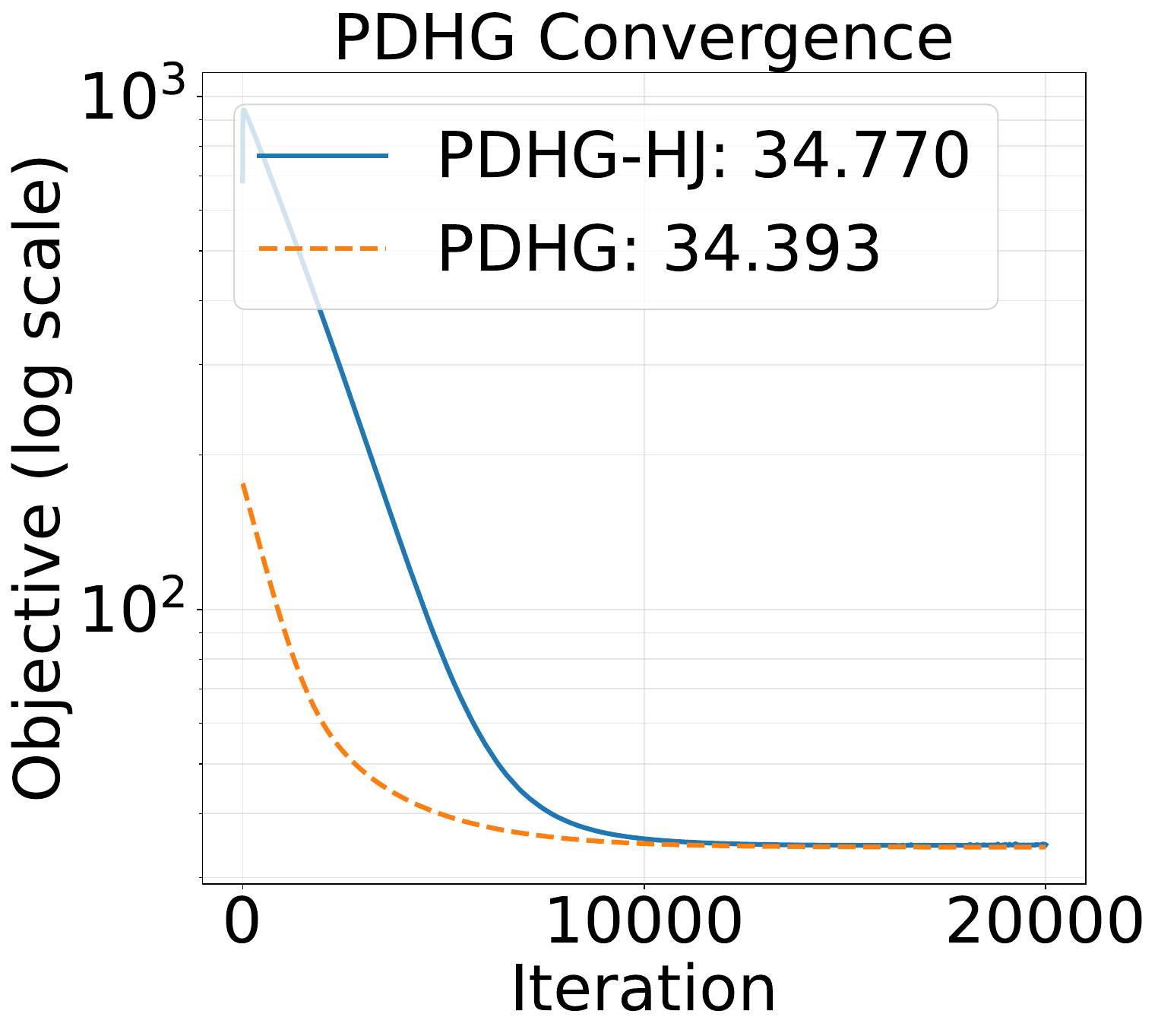}
    \end{tabular}%
    }
    \caption{Total Variation Results}\label{fig: pdhgResults}
\end{figure}

Across all experiments, HJ-Prox tracks the analytical baselines closely and converges to visually indistinguishable solutions. For LASSO and sparse group LASSO, the method performs effective variable selection, shrinking true zero coefficients toward zero as seen in Figures \ref{fig: LASSO and MLR} and \ref{fig: Fused LASSO and Sparse Group LASSO}. In multitask learning and sparse group LASSO, the HJ-Prox iterates closely match the analytical updates. For fused LASSO in Figure \ref{fig: Fused LASSO and Sparse Group LASSO}, HJ-Prox exhibits faster initial convergence but settles farther away than the analytical solver, likely due to differences in formulation (primal vs dual) and the aggressive delta schedule for this particularly challenging proximal operator. We note that convergence speed is problem dependent, as seen in Figures \ref{fig: Fused LASSO and Sparse Group LASSO} and \ref{fig: pdhgResults}, TV and fused LASSO typically require more iterations due to additional subproblems and higher per iteration cost. We note that our goal is not to outperform specialized solvers but to demonstrate that a universal zeroth-order, sampling based proximal approximation integrated with standard splitting algorithms recovers the same solutions with analytical counterparts.

\section{Limitations and Future Work}
Several limitations and promising directions emerge from our analysis. First, our theoretical framework does not account for Monte Carlo sampling errors and assumes exact integral evaluation. To provide more realistic performance guarantees in practical applications, we plan to extend our analysis to incorporate finite-sample errors. Second, convergence can be slow when using set $\delta_k$ schedules and fixed step and sample sizes. Observed in our fused LASSO experiments, adaptive algorithm parameters are often necessary for efficient convergence. We suspect that jointly adapting both the sample size $N$ and $\delta$ throughout the iterative process could be more effective than our current approach of using predetermined schedules as is done for proximal point in~\cite{heaton2024global, ZhangEtAl2024}. We aim to develop adaptive splitting algorithms that dynamically adjust parameters based on problem-specific behavior.
Future work also includes integrating HJ-Prox-based algorithms within a Learning-to-Optimize framework~\cite{chen2022learning, heaton2023explainable, mckenzie2024three, mckenziedifferentiating} to enable automatic tuning of $N$ and $\delta$.
\section{Conclusion}
Our work demonstrates that HJ-Prox can be successfully integrated into operator splitting frameworks while maintaining theoretical convergence guarantees, providing a generalizable method for solving composite convex optimization problems. By replacing exact proximal operators with a zeroth-order Monte Carlo approximation, we have established that algorithms such as PGD, DRS, DYS, and the PDHG method retain their convergence properties under mild conditions. This framework offers practitioners a universal approach to solving complex non-smooth optimization problems, reducing reliance on expensive and complex proximal computations. Our code for experimentation will be available upon publication.  

\newpage

\bibliography{bibliography}
\appendix
\newpage
\section{Proof of HJ-Prox Error Bound}
For completeness and ease of presentation, we restate the theorem.

\textit{\hjproxerrorTheorem{app}}
\begin{proof}
Fix the parameters $t$ and $\delta$. For notational convenience, denote $\prox_{tf}(x)$ by $z^\star(x)$ and $\hjprox_{tf}(x)$ by $z_\delta(x)$. 
Making the change of variable $w = z - z^\star(x)$ in \eqref{eq:hj_prox} enables us to expresses the approximation error as
\begin{eqnarray}
\label{eq:approximation_error}
    z_\delta(x) - z^\star(x) & = & \frac{\int w \exp\left(-\frac{\phi_x(z^\star(x) + w) - \phi_x(z^\star(x))}{\delta}\right) dw}{\int \exp\left(-\frac{\phi_x(z^\star(x) + w) - \phi_x(z^\star(x))}{\delta}\right) dw}.
\end{eqnarray}

Let 
\begin{eqnarray}
    Z_\delta & = & \int \exp\left(-\frac{\phi_x(z^\star(x) + w) - \phi_x(z^\star(x))}{\delta}\right) dw
\end{eqnarray}
and
\begin{eqnarray}
\label{eq:definition_g}
    g(w) & = & \phi_x(z^\star(x) + w ) - \phi_x(z^\star(x)).
\end{eqnarray}
Then
\begin{eqnarray}\label{eq:density}
    \rho_\delta(w) & = & \frac{e^{-\frac{g(w)}{\delta}}}{Z_\delta}
\end{eqnarray} 
defines a proper density. 

Equations \eqref{eq:approximation_error} and \eqref{eq:density} together imply that the approximation error can be written as the expected value of a continuous random variable $W$ whose probability law has the density $\rho_\delta$.
\begin{eqnarray}
\label{eq:error_as_expectation}
    z_\delta(x) - z^\star(x) & = & \int w \, \rho_\delta(w)dw \amp = \amp \E_{\rho_\delta}(W).
\end{eqnarray}
Taking the norm of both sides of \eqref{eq:error_as_expectation} leads to a bound on the norm of the approximation error.
\begin{eqnarray}
\label{eq:size_of_error_bound}
    \|z_\delta(x) - z^\star(x)\| & = & \|\E_{\rho_\delta}(W)\| \amp \leq \amp \E_{\rho_\delta}(\|W\|) \amp \leq \amp \sqrt{\E_{\rho_\delta}\left(\|W\|^2\right)}.
\end{eqnarray}
The first inequality is due to Jensen's inequality since norms are convex. The second is due to the Cauchy-Schwarz inequality. 

Our goal is to show that
\begin{eqnarray}
\label{eq:second_moment_bound}
\sqrt{\E_{\rho_\delta}\left (\lVert W \rVert^2\right)} & \leq & \sqrt{2nt \delta},
\end{eqnarray}
since inequalities \eqref{eq:size_of_error_bound} and \eqref{eq:second_moment_bound} together imply that
\begin{eqnarray}
\|z_\delta(x) - z^\star(x)\| & \leq & \sqrt{2nt \delta}.
\end{eqnarray}

We prove \eqref{eq:second_moment_bound} in two steps. We first show that
\begin{eqnarray}
\label{eq:step_1_bound}
\E_{\rho_\delta}\left( \rVert W \rVert^2\right) & \leq & 2t\E_{\rho_\delta}\left(\left\langle W, \nabla g(W)\right\rangle \right),
\end{eqnarray}
where $g$ is the convex function defined in \eqref{eq:definition_g}.
We then show that
\begin{eqnarray}
\label{eq:step_2_equality}
\E_{\rho_\delta}\left(\left\langle W, \nabla g(W)\right\rangle \right) & = & n\delta.
\end{eqnarray}
Before proceeding to prove these steps, we address our abuse of notation in \eqref{eq:step_1_bound} and \eqref{eq:step_2_equality}. 
Although $\nabla g$ may not exist everywhere, it exists almost everywhere. 
Recall that $g$ is locally Lipschitz because it is convex. Furthermore, any locally Lipschitz function is differentiable almost everywhere by Rademacher's theorem. Hence $\nabla g$ exists almost everywhere. Consequently, the expectation $\E_{\rho_\delta}\left \langle W, \nabla g(W) \right\rangle$ is well defined.

To show \eqref{eq:step_1_bound}, first note that by Fermat's rule, $0 \in \partial\phi_x(z^\star(x))$ where $\partial \phi(z)$ denotes the subdifferential of $\phi$ at $z$.
Consequently, for any $z \in \Real^n$
\begin{equation}
\label{eq:Inequality}
\begin{split}
    \phi_x(z) & \amp \geq \amp  \phi_x(z^\star(x)) + \langle 0, z - z^\star(x)\rangle + \frac{1}{2t}\|z-z^\star(x)\|^2 \\
    & \amp = \amp \phi_x(z^\star(x)) + \frac{1}{2t}\|z-z^\star(x)\|^2,
\end{split}
\end{equation}
since $\phi_x(z)$ is $\frac{1}{t}$-strongly convex. Plugging $z = z^\star(x) + w$ into inequality \eqref{eq:Inequality} implies that
\begin{eqnarray}
\label{eq:strong_convexity_inequality}
g(w) & \geq & \frac{1}{2t}\lVert w \rVert^2.
\end{eqnarray}
Suppose $\nabla g$ exists at a point $w$. Then
\begin{eqnarray}
g(0) & \geq & g(w) + \langle \nabla \ g(w), 0 - w \rangle,
\end{eqnarray}
because $g$ is convex. Since $g$ vanishes at zero, rearranging the above inequality gives
\begin{eqnarray}
\label{eq:intermediate_inquality}
\langle \nabla g(w), w \rangle & \geq & g(w).
\end{eqnarray}
Inequalities \eqref{eq:strong_convexity_inequality} and \eqref{eq:intermediate_inquality} give
\begin{eqnarray}
\frac{1}{2t}\lVert w \rVert^2 &  \leq  & \langle \nabla g(w), w \rangle,
\end{eqnarray}
which implies that
\begin{eqnarray}
    \frac{1}{2t}\int \lVert w \rVert^2 \rho_\delta(w)dw & \leq & \int \langle w,\nabla g(w)\rangle \rho_{\delta}(w)dw.
\end{eqnarray}

To show \eqref{eq:step_2_equality}, consider $w$ where $\nabla g$ exists and let $h(w) = \exp(-g(w)/\delta)$. By the chain rule
\begin{eqnarray}
\label{eq:chain_rule}
\frac{\partial}{\partial w_j} g(w)h(w) & = & - \delta \frac{\partial}{\partial w_j}h(w).
\end{eqnarray}
Integrating both sides of \eqref{eq:chain_rule} over $\Real^d$ gives
\begin{equation}
\label{eq:component}
\begin{split}
    \int_{\Real^n} w_j \frac{\partial}{\partial w_j}g(w) h(w) dw  & \amp = \amp  -\delta \int_{\Real^n}  w_j \frac{\partial}{\partial w_j} h(w) dw \\ 
    & \amp = \amp -\delta\int_{\Real^{n-1}}\left[\int_{-\infty}^\infty w_j \frac{\partial}{\partial w_j}h(w)dw_j\right]dw_{-j},
\end{split}
\end{equation}
where
$w_{-j}$ is the subvector of $w$ containing all but its $j$th element.

Applying integration by parts on the right hand side of \eqref{eq:component} gives
\begin{eqnarray}
\int_{-\infty}^\infty w_j \frac{\partial}{\partial_{w_j}}h(w) d w_j & = & w_j h(w)\Big |_{-\infty}^\infty - \int_{-\infty}^\infty h(w)dw_j.
\end{eqnarray}
Note that \eqref{eq:strong_convexity_inequality} implies that
\begin{eqnarray}
\underset{w_j \rightarrow \infty}{\lim}\;
\left\lvert w_j h(w) \right\rvert & \leq & \underset{w_j \rightarrow \infty}{\lim}\;
\left\lvert w_j  e^{-\frac{\lVert w \rVert^2}{2t}}\right\rvert \amp = \amp 0.
\end{eqnarray}
Consequently,
\begin{eqnarray}
\label{eq:integration_by_parts}
       \int_{\Real^{d}} w_j \frac{\partial}{\partial_{w_j}}h(w)dw & = & \int_{\Real^{d-1}}\left[-\int_{-\infty}^\infty h(w)dw_j \right]dw_{-j} \amp = \amp -Z_\delta.
\end{eqnarray}
Equations \eqref{eq:component} and \eqref{eq:integration_by_parts} together imply that
\begin{eqnarray}
\E_{\rho_\delta}\left( W_j \frac{\partial}{\partial w_j} g(W) \right) & = & \delta.
\end{eqnarray}
The linearity of expectations gives \eqref{eq:step_2_equality} completing the proof.
\end{proof}

\section{HJ-Prox-based PGD Convergence}
For completeness and ease of presentation, we restate the theorem.

\noindent \textbf{Proof of Thm.~\ref{theorem:hjpgd}.}
\textit{\pgdTheorem{app}}
\begin{proof}
For appropriately chosen step-size $t$, the PGD algorithm map is averaged and its fixed points coincide with the global minimizers of $f$ (as shown in the Lemma below).
\begin{lemma}[Averagedness and Fixed Points of PGD]\label{lem: PGD fixed-point}
Let $0<t<\frac{2}{L}$ and define, for $x \in \Real^n$
\begin{eqnarray}
T(x) & = & \prox_{t g}\bigl(x - t\nabla f(x)\bigr)
.
\end{eqnarray}
Then $T$ is an averaged operator, and its fixed points $ \text{Fix} (T)$ coincide with $f$'s global minimizers $X^*$ \cite[Section 4.2]{ParikhBoyd2014}.
\end{lemma}
The PGD iterates are computed by applying the mapping 
$T(x) = \prox_{t g}(x - t \nabla f(x))$.
By Lemma~\ref{lem: PGD fixed-point}, $T$ is an averaged operator and $x_k \rightarrow x^* \in \text{Fix} (T) =  \mathnormal{X}^*$. 

The HJ-PGD iterates can be written as 
\begin{eqnarray}
    \hat{x}_{k+1} & = & \prox_{t g}^{\delta_k}(\hat{x}_k - t \nabla f(\hat{x}_k)) \amp = \amp T(\hat{x}_k) + \varepsilon_k,
\end{eqnarray}
where
\begin{eqnarray}
\varepsilon_k & = & \prox_{t g}^{\delta_k}(\hat{x}_k - t \nabla f(\hat{x}_k)) - \prox_{t g }(\hat{x}_k - t \nabla f(\hat{x}_k) ).
\end{eqnarray}
Since $\sum_{k} \sqrt{\delta_k}$ is finite, $\sum_k^\infty \|\varepsilon_k\|$ is finite by Theorem~\ref{thm:hj_prox_error_bound}. Consequently, $\hat{x}_k \to x^\star\in X^*$ by Theorem~\ref{theorem:perturbed_KM}.
\end{proof}

\section{HJ-Prox-based DRS Convergence}
We restate the statement of the theorem for readability.

\noindent \textbf{Proof of Thm.~\ref{theorem:hjdrs}.}
\textit{\drsTheorem{app}}
\begin{proof}
The DRS algorithm map is averaged and its fixed points coincide with the global minimizers of $f$.

\begin{lemma}[Averagedness and Fixed Points of DRS]\label{lem: DR-fixed-point}
Let $t>0$ and define, for $z\in\Real^n$
\begin{eqnarray}
\label{eq:DRS_update}
T(z) & = & z + \prox_{t g}\!\bigl(2\,\prox_{t f}(z) - z\bigr) - \prox_{t f}(z).
\end{eqnarray}
Note this is the fixed point operator for the dual variable in the DRS algorithm.
Then $T$ is firmly nonexpansive (hence averaged), and
\begin{eqnarray}
\text{Fix}(T) & = & \{ z: \prox_{t f}(z) \in \mathnormal{Z}^* \,\}. 
\end{eqnarray} \cite[Remark 5]{LionsMercier1979}
\end{lemma}
By Lemma~\ref{lem: DR-fixed-point}, $z_k \rightarrow z^*$ and $\prox(z^*) = x^* \in \mathnormal{X}^*$. 
We can express the HJ-DRS update in terms of the DRS algorithm map $T$ \eqref{eq:DRS_update}.
\begin{eqnarray}
    \hat{z}_{k+1}
    &=& T(\hat{z}_k) + \varepsilon_k,
\end{eqnarray}
where
\begin{eqnarray}
\varepsilon_k & = & \prox_{t h}(w_k + 2\kappa_k) - \prox_{t h}(w_k) +\zeta_k -\kappa_k,
\\
w_k & = & 2\prox_{t g}(\hat{z}_k) - \hat{z}_k,
\end{eqnarray}
and
\begin{eqnarray}
    \zeta_k & = & \prox^{\delta_k}_{t g}(\hat{z}_k) - \prox_{t g}(\hat{z}_k) \\
    \kappa_k & = & \prox^{\delta_k}_{t f}(\hat{z}_k) - \prox_{t f}(\hat{z}_k).
\end{eqnarray}
We have the following bound
\begin{eqnarray}
    \|\varepsilon_k\| & \leq & \|w_k+ 2\kappa_k -w_k\| + \|\zeta_k\| +\|\kappa_k\| \amp = \amp  3\|\kappa_k\| + \|\zeta_k\|,
\end{eqnarray}
which follows from the triangle inequality and the fact that proximal mappings are nonexpansive.

Since $\sum_{k} \sqrt{\delta_k}$ is finite $\sum_k^\infty \|\varepsilon_k\|$ is finite by Theorem~\ref{thm:hj_prox_error_bound}. Consequently, $\hat{z}_k \to z^\star$ by Theorem~\ref{theorem:perturbed_KM}. Since proximal maps are continuous, $\hat{x}_k = \prox(\hat{z}_k) \rightarrow \prox(z^\star) \in X^*$.
\end{proof}
\section{HJ-Prox-based DYS Convergence}
For completeness and ease of presentation, we restate the theorem.

\noindent \textbf{Proof of Thm.~\ref{theorem:hjdys}.}
\textit{\dysTheorem{app}}
\begin{proof}
For appropriately chosen step-size $t$, the DYS algorithm map is averaged and its fixed points coincide with the global minimizes of $f+g+h$.
\begin{lemma}[Averagedness and Fixed Points of DYS]\label{lem: DYS-fixed-point}
Let $t > 0$ and define, for $z \in \Real^n$,
\begin{eqnarray}
\label{eq:DYS_update}
    T(z) = z - \prox_{t f}(z) + \prox_{t g}\bigl(2\prox_{t f}(z)-z-t\nabla h(\prox_{t f}(z)\bigr).
\end{eqnarray}
Note this is the fixed point operator for the DYS algorithm and its fixed points Fix(T) coincide with global minimizers $X^*$. T is firmly nonexpansive (hence averaged), and 
\begin{eqnarray}
    Fix(T) = \{z : x \in X^\star\}
\end{eqnarray}\cite[Theorem 3.1]{DavisYin2015}
\end{lemma}
By Lemma ~\ref{lem: DYS-fixed-point}, $z_k \rightarrow z^\star$ and $z^* \in \mathnormal{X}^*$. We can express the HJ-DYS update in terms of DYS algorithm map $T$ \eqref{eq:DYS_update}.
\begin{eqnarray}
    \hat{z}_{k+1} &=& T(\hat{z}_k) + \varepsilon_k,
\end{eqnarray}
where 
\begin{eqnarray}
    \varepsilon_k &=& \prox_{tg}\big({S}_t(z_k) + d_k\big) - \prox_{t g}\big({S}_t(z_k)\big) + \zeta_k - \kappa_k\\
    {S}_t(z_k) &=& 2\prox_{t f}(z_k) - z_k - t \nabla h\bigl(\prox_{tf}(z_k)\bigr)\\
    d_k &=& 2 \kappa_k - t [\nabla h \bigl(\prox_{tf}(z_k) +\kappa_k\bigr)- \nabla h \bigl(\prox_{t g}(z_k)\bigr)]
\end{eqnarray}
and 
\begin{eqnarray}
    \zeta_k &=& \prox^{\delta_k}_{t g}(\hat{z}_k) - \prox_{t g}(\hat{z}_k) \\
    \kappa_k &=& \prox^{\delta_k}_{t f}(\hat{z}_k) - \prox_{t f}(\hat{z}_k).
\end{eqnarray}
We have the following bound
\begin{eqnarray}
    \|\varepsilon_k\|  & \leq &  (1+ t L)\|\kappa_k\| + \|\zeta_k\|, 
\end{eqnarray}
which follows from the triangle inequality, $L$-smoothness of $h$, and the fact that proximal mappings are nonexpansive. 

Since $\sum_{k} \sqrt{\delta_k}$ is finite $\sum_k^\infty \|\varepsilon_k\|$ is finite by Theorem~\ref{thm:hj_prox_error_bound}. Consequently, $\hat{z}_k \to z^\star$ where $z^\star$ is a global minimizer of $f+g+h$ by Theorem~\ref{theorem:perturbed_KM} and Lemma ~\ref{lem: DYS-fixed-point}.
\end{proof} 
\section{HJ-Prox-based PDHG Convergence}
For completeness and ease of presentation, we restate the theorem.

\noindent \textbf{Proof of Thm.~\ref{theorem:hjpdhg}.}
\textit{\pdhgTheorem{app}}
\begin{proof}
For appropriately chosen $\tau,\sigma$ the PDHG algorithm map is averaged and its fixed points corresponding to $x_k$ updates coincide with the global minimizes of $f(x)+g(Ax)$.
\begin{lemma}[Averagedness and Fixed Points of PDHG]\label{lem: PDHG-fixed-point}
Let $\tau,\sigma > 0$ satisfying $\tau \sigma \|A\|^2 < 1$ and define, for $z \in \Real^m$ and $w \in \Real^n$ 
\begin{eqnarray}
\label{eq:PDHG_update}
    T(z,w) = \begin{bmatrix}
\prox_{\tau f}\!\big(z - \tau A^\top \prox_{\sigma g^*}(w + \sigma Az)\big) \\
\prox_{\sigma g^*}(w + \sigma Az)
\end{bmatrix}.
\end{eqnarray}
Let $V = \mathrm{diag}(\tfrac{1}{\tau}I_n,\, \tfrac{1}{\sigma}I_m)$. On a product space with a weighted inner product
$\langle (x,y),(x',y')\rangle_V = \tfrac{1}{\tau}\langle x,x'\rangle + \tfrac{1}{\sigma}\langle y,y'\rangle$,
the map T is an averaged operator. 
Note this is the fixed point operator for the PDHS algorithm and its fixed points Fix(T) coincide with the set of primal-dual KKT saddle points for $f(x)+g(Ax)$, where the primal point coincides with the global minimizers $X^*$. T is firmly nonexpansive (hence averaged), and 
\begin{eqnarray}
    Fix(T) = \{(z^*,w^*) : z^* \in X^\star\}.
\end{eqnarray}\cite[Algorithm 1, Thm. 1]{ChambollePock2011} \cite[Lemma 2]{Fercoq2022}
\end{lemma}
By Lemma ~\ref{lem: PDHG-fixed-point}, $z_k \rightarrow z^\star$ and $z^* \in \mathnormal{X}^*$. We can express the HJ-PDHG update in terms of PDHG algorithm map $T$ \eqref{eq:PDHG_update}.
\begin{eqnarray}
    (\hat{z}_{k+1},\hat{w}_{k+1}) &=& T(\hat{z}_{k},\hat{w}_{k}) + \varepsilon_k
\end{eqnarray}
where 
\begin{eqnarray}
    \varepsilon_k &=& \begin{bmatrix}
\prox_{\tau f}(u_k - \tau A^\top \zeta_k) - \prox_{\tau f}(u_k)+\kappa_k \\
\zeta_k
\end{bmatrix}\\
u_k &=& \hat{z}_k - \tau A^\top \prox_{\sigma g^\star}(w_k+ \sigma A\hat{z}_k),
\end{eqnarray}
and 
\begin{eqnarray}
    \zeta_k &=& \prox^{\delta_k}_{\sigma g^*}(\hat{z}_k + \sigma A\hat{w}_k) - \prox_{\sigma g^*}(\hat{z}_k + \sigma A\hat{w}_k) \\
    \kappa_k &=& \prox^{\delta_k}_{\tau f}(\hat{w}_k - \tau A^\top \hat{z}_k) - \prox_{\tau f}(\hat{w}_k - \tau A^\top \hat{z}_k)
\end{eqnarray}
In the weighted norm $\|(w,z)\|_V^2 = \frac{1}{\tau} \|w\|^2 + \frac{1}{\sigma}\|z\|^2$ , we have the following bound
\begin{eqnarray}
    \| \varepsilon_k\|^2_V &\leq& \big(2\tau\|A\|_{\text{op}}^2 + \frac{1}{\sigma}\big)\|\zeta_k\|^2 + \frac{2}{\tau}\|\kappa_k\|^2
\end{eqnarray}
which follows from the fact that proximal mappings are nonexpansive and from $\|A^\top\|_{\text{op}} = \|A\|_{\text{op}}$.
Since $\sum_{k} \sqrt{\delta_k}$ is finite $\sum_k^\infty \|\varepsilon_k\|$ is finite by Theorem~\ref{thm:hj_prox_error_bound}. Consequently, $(z_k,w_k) \rightarrow (z^*,w^*)$ where $z^\star$ is a global minimizer of $f+g$ by Theorem~\ref{theorem:perturbed_KM} and Lemma ~\ref{lem: PDHG-fixed-point}.
\end{proof}

\section{Experiment Details}\label{AP: Experiments}
HJ-Prox and analytical counterparts run through all iterations. Every experiment simulates a ground truth structure with added noise and blur depending on problem setup. All parameters and step sizes are matched between HJ-Prox and the analytical counterparts to ensure a fair comparison. The HJ-Prox $\delta$ sequence follows a schedule
\begin{eqnarray}
    \delta_k & = & \frac{1}{k^{2.00001}},
\end{eqnarray}
where $k$ denotes iteration number. The defined schedule decays strictly faster than $1/k^2$ satisfying conditions used in Theorem~\ref{theorem:perturbed_KM}.
\subsection{PGD: LASSO Regression}
We solve the classic LASSO regression problem using PGD. The simulation setup involves a design matrix $X \in \mathbb{R}^{250 \times 500}$ with 250 observations and 500 predictors. The true coefficients $\beta$ are set such that $\beta^{400:410} = 1$ and all others are zero. The objective function is written as,
\begin{equation}
    \underset{\beta}{\arg\min}\ \frac{1}{2}\|X\beta - y\|_2^2 + \lambda\|\beta\|_1
\end{equation}
\[
X\in\mathbb{R}^{250\times 500},\quad
\beta\in\mathbb{R}^{500},\quad
y\in\mathbb{R}^{250}.
\]
The analytical PGD baseline performs a gradient step on the least-squares term followed by the exact soft thresholding.
\subsection{DRS: Multitask Regression}
Multitask regression learns predictive models for multiple related response variables by sharing information across tasks to enhance performance. We solve this problem using Douglas-Rachford splitting, employing HJ-Prox in place of analytical updates.
We group the quadratic loss with the nuclear norm regularizer to form one function and the row and column group LASSO terms to form the other. Both resulting functions are non-smooth, requiring HJ-Prox for their proximal mappings. The simulation setup involves $n=50$ observations, $p=30$ predictors, and $q=9$ tasks. The objective function is written as, 
\begin{equation}
\underset{B}{\arg\min}\;
\frac{1}{2}\,\| X B-Y\|_{F}^{2}
\;+\;\lambda_{1}\,\| B\|_{*}
\;+\;\lambda_{2}\sum_{i}\| b_{i,\cdot}\|_{2}
\;+\;\lambda_{3}\sum_{j}\|b_{\cdot,j}\|_{2}
\end{equation}
\[
X\in\mathbb{R}^{50\times 30},\quad
B\in\mathbb{R}^{30\times 9},\quad
Y\in\mathbb{R}^{50\times 9}.
\]
The analytical counterpart for Douglas Rachford Splitting utilizes singular value soft thresholding for the nuclear norm and group soft thresholding for the row and column penalties. These regularizers are integrated with fast iterative soft thresholding (FISTA) to handle the data fidelity term with nuclear norm regularization and Dykstra's algorithm to handle the sum of row and column group LASSO penalties.

\subsection{DRS: Fused LASSO}
The fused LASSO is commonly used in signal processing to promote piecewise smoothness in the solution. We apply it to recover a Doppler signal with length $n=256$ using a third-order differencing matrix $D$. We solve this problem with DRS, comparing two implementation strategies: an exact method using product-space reformulation motivated by \cite{TibshiraniTaylor2011}, and an approximate method using HJ-Prox. The objective function is written as, 

\begin{equation}
    \underset{B}{\arg\min}\ \frac{1}{2}\|\beta - y\|^2 + \lambda\|D\beta\|_1
\end{equation}
\[
\beta \in\mathbb{R}^{256},\quad
y\in\mathbb{R}^{256},\quad
D \in\mathbb{R}^{253\times 256}.
\]

The proximal operator of $\lambda\|D\beta\|_1$ has no closed-form solution for general linear operators $D$. The analytical counterpart addresses this by reformulating the problem in a product space with auxiliary variable $w = D\beta$, yielding separable proximal operators (weighted averaging and soft thresholding) at the cost of inverting terms including $D^\top D$ at each iteration. In contrast, our HJ-Prox variant directly approximates the intractable proximal operator through Monte Carlo sampling. As a reminder, both use identical DRS parameters for fair comparison. 

\subsection{DYS: Sparse Group LASSO}
The sparse group LASSO promotes group-level sparsity while allowing individual variable selection within groups, which is useful when certain groups are relevant but contain unnecessary variables.
We solve this problem using DYS, employing HJ-Prox for the proximal operators of the non-smooth regularizers. The simulation setup involves $n=300$ observations with $G=6$ groups, each having 10 predictors. The objective function is written as, 
\begin{equation}
\underset{\beta}{\arg\min}\;
\frac{1}{2}\,\| X \beta - y\|_{2}^{2}
\;+\;\lambda_1\sum_{g=1}^{G}\| \beta_{g}\|_{2}
\;+\;\lambda_2\| \beta\|_{1}
\end{equation}
\[
X \in \mathbb{R}^{300 \times 60},\quad
\beta \in\mathbb{R}^{60},\quad
y\in\mathbb{R}^{300}.
\]

The analytical counterpart for DYS solves the sparse group LASSO by using soft thresholding for the $\ell_1$ penalty and group soft thresholding for the group $\ell_2$ penalty.
\subsection{PDHG: Total Variation}
Lastly, we implement PDHG method to solve the isotropic total variation regularized least‐squares problem. We apply the proximal operator for the data fidelity term via its closed‐form update and employ our HJ‐based proximal operator for the total variation penalty. For this experiment, we recover a smoothed 64 x 64 black and white image from a noisy and blurred image $y$. The objective function is written as, 

\begin{equation}
\underset{\beta}{\arg\min}\;
\frac{1}{2}\,\| X \beta - y\|_{F}^{2}
+\lambda \text{TV}(\beta)
\end{equation}
\[
\beta \in\mathbb{R}^{64 \times 64},\quad
y\in\mathbb{R}^{64 \times 64}.
\]

The (slightly smoothed) isotropic TV we use to evaluate the objective is
\begin{eqnarray}
    \mathrm{TV}(\beta)
 & = &
\sum_{i=1}^{64}\sum_{j=1}^{64}
\sqrt{\,(\nabla_x \beta)_{i,j}^2 + (\nabla_y \beta)_{i,j}^2}.
\end{eqnarray}

The analytical counterpart for PDHG  algorithm updates dual variables using closed-form scaling for data fidelity and clamping (for $\ell_2$ projection of TV dual), and primal variables using Fast Fourier transform convolution and divergence via finite differences.

\end{document}